\documentclass[reqno,11pt]{amsart}

\newtheorem{theorem}{Theorem}[section]
\newtheorem{lemma}[theorem]{Lemma}
\newtheorem{corollary}[theorem]{Corollary}

 \theoremstyle{definition}

\newtheorem{assumption}{Assumption}[section]
\newtheorem{example}[theorem]{Example}

\theoremstyle{remark}
\newtheorem{remark}[theorem]{Remark}

\newcommand\bR{\mathbb{R}}

\newcommand\bZ{\mathbb{Z}}

\newcommand\cB{\mathcal{B}}

\newcommand\cF{\mathcal{F}}
\newcommand\cG{\mathcal{G}}

\newcommand\cP{\mathcal{P}}

\newcommand\frN{\mathfrak{N}}

\newcommand{\WO}{\overset{\scriptscriptstyle0}%
{W}\,\!}

\newcommand{\tr}{\text{\rm tr}\,}

\renewcommand\){{\rm)}}

\newcommand{\mysection}[1]{\section{#1}
\setcounter{equation}{0}}

\newcommand\cbrk{\text{$]$\kern-.15em$]$}}
\newcommand\opar{
\text{\,\raise.2ex\hbox{${\scriptstyle |}$}\kern-.34em$($}}
\newcommand\cpar{%
\text{$)$\kern-.34em\raise.2ex\hbox{${\scriptstyle |}$}}\,}
\newcommand\obrk{\text{$[$\kern-.15em$[$}}

\newcommand\R{\mathbb R}
\begin{document}
 
\title[Parabolic estimates and Poisson process]
{Poisson stochastic process and basic Schauder and
Sobolev  estimates  in the theory of parabolic
equations}
\author{N.V. Krylov}%
\thanks{The first author was partially supported by
 NSF Grant DMS-1160569  
and by a grant 
from the Simons Foundation (\#330456 to Nicolai Krylov)}
\address{127 Vincent Hall, University of Minnesota, Minneapolis,
  MN, 55455}
\email{nkrylov@umn.edu}

\author{E. Priola}
\thanks{   The second author
was partially   supported by the Italian PRIN project 2010MXMAJR}
\address{Dipartimento di Matematica, 
  Universit\`a  di Torino,
 Via Carlo Alberto
10, 10123 Torino, Italy} \email{enrico.priola@unito.it}


\begin{abstract}
We show among other things how  knowing Schauder  or
Sobolev-space estimates for the one-dimensional
heat equation allows one to derive their
multidimensional analogs for equations with
coefficients depending only on time variable
with the {\em same\/}
constants as in the case of the one-dimensional heat equation.
The method is   quite general and is  based on using the Poisson
stochastic process.  It also applies to equations involving
non-local operators.  It looks like no other method is
available at this time and it is a very challenging problem
to find a purely analytic approach to proving such results. 
 
\end{abstract}

\maketitle

\mysection{Introduction}
                                       \label{section 6.29.1}

In this paper we present a method allowing one, in particular, 
to obtain various estimates for the multidimensional
second-order
parabolic equations of main type with time dependent coefficients
with the {\em same\/} constants
as in the case of the one-dimensional heat equation,
provided that the matrix of the second-order coefficients
dominates the  identity  matrix.

The method is universal in the sense that it works  in the same
way
for H\"older- or Sobolev-space estimates,  for scalar equations and   even 
for not necessarily parabolic systems. 
The main condition for it to work is that the equations
should be commuting with space translations 
 (more generally, should be commuting with a 
commutative group 
of affine mappings)  and the estimates
should be space-translation invariant as well.

We start with Section \ref{section 1.25.1}
and show our main idea on the example of deriving
basic Schauder and Sobolev-space estimates
for the heat equation in 2 space dimension
from the similar estimates for the 
heat equation in 1 space dimension.
Here we just use the Poisson process.

In Section \ref{section 6.29.2} we show how the method
works for multidimensional parabolic equations
with measurable coefficients
depending only on the time variable, provided
that the matrix of the coefficients dominates the 
 identity  matrix.
This time an integral of nonrandom functions against
the Poisson process is involved.

As a corollary we obtain that for {\em elliptic\/} equations of main type
with constant coefficients the constant in the estimate 
of the $C^{\alpha}$-semi-norm of  the second-order {\em directional\/}
derivatives of solutions through 
the $C^{\alpha}$-semi-norm of the free term
is independent of the space dimension. The same is also noted
for the $L_{p}$-estimates of the second-order {\em directional\/}
derivatives of solutions through the $L_{p}$-norm 
of the free term.

In Section \ref{section 6.29.3} we present our method
in a more abstract form for evolution equations
when  the norms
are not necessarily translation invariant, but invariant
relative to a group of affine mappings of the space
and the equations commute with that group.
In Example \ref{example 2.29.4} we show a result of application
of our general theorem, Theorem \ref{theorem 2.29.1},
   which allows us to obtain the
Schauder estimates for a parabolic equation with
{\em space-dependent\/} coefficients with the same
constants as in the case of the 2 dimensional heat equation.
In Example \ref{example 5.6.1} we apply Theorem \ref{theorem 2.29.1}
  to a hyperbolic system. 
In Example \ref{example 10.6.1} we show an application
of our results to the hyperbolic systems from \S 7.3.3
of Evans's book \cite{Ev}.

 Section \ref{section 6.29.5} contains the proof
of Theorem \ref{theorem 2.27.1}, which is used
in Section \ref{section 6.29.6} to prove Theorem \ref{theorem 2.29.1}.
  Finally, in Section 7 we present an 
extension of our method to treat non-local operators.
 
The origin of our ideas lies in the theory
of stochastic partial differential equations (SPDEs)
and can be found in 
the proof of Theorem 2.1 of \cite{Kr_96}.
This idea can be implemented quite formally without
using the theory of SPDEs, see, for instance,
\cite{K1}
and \cite{P}, where still one needs to be familiar
with the It\^o stochastic integral with respect to the Wiener
process
 albeit of nonrandom
functions. 

It turns  out  that replacing the Wiener
process with the Poisson process in the original idea
leads to much simpler SPDEs which, actually,
are just usual equations
with  discontinuities in time at random well separated
moments, dealing with which does
not require any knowledge of stochastic
integration. Turning to  the  Poisson processes
has also an advantage that we can consider integro-differential
equations   (cf. Theorem \ref{theorem 9.17.2}).  

At the same time we can easily recover the results
obtained by using methods in \cite{K1}
and \cite{P}. The probabilistic reason (which is not used
in the article) for that lies in
 the well-known central limit
theorem according to which $(2\lambda)^{-1/2}(\pi^{\lambda,1}_{t}
-\pi^{\lambda,2}_{t})$ tends in law to $w_{t}$
as $\lambda\downarrow0$, where
 $\pi^{\lambda,i}_{t}$, $i=1,2$, are independent Poisson
processes with intensity $\lambda$  and $w_{t}$ is a Wiener
process.

In conclusion we note that the scope of applications
of  Theorems 
\ref{theorem 2.27.1} and
 \ref{theorem 2.29.1} is much wider than only
the examples given in the article. For instance,
one could consider integro-differential equations 
or higher order equations, or else the combinations of those.
We plan to explore these possibilities in the
near future.

In the whole article $T$ is a fixed number in $(0,\infty)$,
$\bR^{d}$ is a  Euclidean space of points $x=(x^{1},...,x^{d})$,
$x^{1},...,x^{d}\in(-\infty,\infty)$, 
$S_{1}:=\{x \in \R^d : \, |x|=1\}$ is the unit sphere, and 
the standard stipulation
about the summation with respect to repeated indices is
enforced. Also we use standard notation for derivatives,
spaces, semi-norms, and norms which can be found in
\cite{K3}, \cite{LSU}, \cite{Li}.  We only recall
  what H\"older functions spaces are. 
  The space 
$C^{\alpha}(\bR^d)$, $\alpha \in (0,1)$,  is the 
space of all 
real-valued functions $f$ on $\bR^{d}$ 
for which the following norm
$$
 \| f \|_{C^{\alpha}(\bR^d)} = \sup_{x \in \bR^d} |f(x)| +
  [f]_{C^{\alpha}(\bR^d)}
$$
 is finite, where
$$
  [f]_{C^{\alpha}(\bR^d)} =
   \sup_{x \not = y } \frac{ |f(x) - f(y)|}{ |x-y|^{\alpha}}.
$$
 As usual, by $C^{2+ \alpha}(\bR^d)$
we mean the space of real-valued 
twice continuously differentiable
functions $f$
on $\bR^{d}$ having
 finite norm
$$
\| f \|_{C^{2+\alpha}(\bR^d)} = \sup_{x \in \bR^d}
( |f(x)|+|Df(x)|+|D^{2}f(x)|)+[D^{2}f]_{C^{\alpha}(\bR^d)} ,
$$
where $Df$ is the gradient of $f$ and $D^{2}f$ is its
Hessian.

\mysection{One dimensional heat equation}
                                     \label{section 1.25.1}

Consider the   
problem of solving the equation
\begin{equation}
                                      \label{1.23.5}
\partial_{t}u(t,x)=D^{2}u(t,x)+f(t,x)
\end{equation}
for $t\in(0,T)$, $x\in \bR$, with zero initial condition,   i.e., $u(0,
 \cdot ) =0$.   To be more precise we treat the problem
in the integral form:
\begin{equation}
                                                  \label{ci1}
   u(t,x) = \int_0^t (D^2u(s,x) + f(s,x))ds,\;\; t 
\in [0,T],\; x \in \bR. 
 \end{equation}
 
For a real-valued function $f(t,x)$, $t\in(0,T)$, 
$x\in \bR^{d}$,
write 
$$
f\in B_{c}((0,T),C^{\infty}_{0} (\bR^{d}))
$$
if $f$ is a Borel bounded function, such that
$f(t,\cdot)\in C^{\infty}_{0}(\bR^{d})$ for any $t$,
for any $n=0,1,...$, the
$C^{n} (\bR^{d})$-norms of $f(t,\cdot)$ are bounded   on $(0,T)$,
and the supports of $f(t,\cdot)$ belong to the same ball.

  Fix $\alpha\in(0,1)$ and $p\in(1,\infty)$.
 One knows (see, for instance,
\cite{K3}, \cite{LSU}, \cite{Li})
that if $f\in B_{c}((0,T),C^{\infty}_{0}(\bR))$,
then the above problem has a 
solution  
 $u(t,x)$
having the following properties:

(a) it is continuous in $[0,T]\times\bR$;

(b)
  $u(t, \cdot ) \in C^{2+\alpha}(\bR)$, for any $t \in [0,T]$, and 
\begin{equation}
                                      \label{1.23.1}
  \sup_{ t \in [0,T]  }\|u(t,\cdot)\|_{C^{2+\alpha}(\bR)}
\leq N_{0}(T,\alpha) \sup_{  t \in (0,T)  }
\|f(t,\cdot)\|_{C^{\alpha}(\bR)},
\end{equation}
where $N_{0}(T,\alpha)$ is a (finite)
 constant depending only on $T$ and $\alpha$.
There is only one solution with these properties and,
furthermore,  
\begin{equation}
                                      \label{1.22.2}
\sup_{(t,x) \in [0,T]\times \bR}|u(t,x)|
\leq  T \sup_{(t,x) \in (0,T)\times \bR}|f(t,x)|,
\end{equation} 
\begin{equation}
                                      \label{1.22.1}
\sup_{t \in [0,T]}[D^{2}u(t,\cdot)]_{C^{\alpha}(\bR)}
\leq N_{0}(\alpha) \sup_{t \in (0,T)}[f(t,\cdot)]_{C^{\alpha}(\bR)},
\end{equation}
\begin{equation}
                                      \label{1.24.6}
\|D^{2}u\|^{p}_{L_{p}((0,T)\times\bR)}
\leq N_{p}\|f\|^{p}_{L_{p}((0,T)\times\bR)},
\end{equation}
  where $L_p$-spaces are defined with respect to 
 Lebesgue measure 
and $N_{0}(\alpha),N_{p}$ are some constants.

Take a sequence $\tau_{1}=\tau_{1}(\omega),\tau_{2}
=\tau_{2}(\omega),...$
of independent random variables defined on a probability
space $(\Omega,\cF,P)$ with common 
exponential distribution with parameter $\lambda>0$,
so that $P(\tau_{n}>t)=e^{-\lambda t}$ for $t\geq0$
and $n=1,2....$.
Define
$$
\sigma_{0}=0,\quad
\sigma_{n}=\sum_{i=1}^{ n}\tau_{i},\quad n=1,2,..., 
\quad \pi_{t}=\pi_{t}(\omega)=\sum_{n=1}^{\infty}I_{\sigma_{n}\leq t}
$$
  (where $I_{\sigma_{n}\leq t}$ denotes 
the indicator function of the event 
$\{\sigma_{n}\leq t \}$). 
We see that $\pi_{t}$ is the  number of consecutive
sums of $\tau_{i}$ which lie on $[0,t]$. The 
counting process $\pi_{t}$ is known as a Poisson process
with parameter $\lambda$,
 for $0\leq s\leq t<\infty$ and $k=0,1,...$
it holds that
$$
P(\pi_{t}-\pi_{s}=k)=\frac{[\lambda(t-s)]^{k}}{k!}
e^{-\lambda(t-s)},
$$
and, moreover, $\pi_{t}-\pi_{s}$ is independent 
of the trajectory $\{\pi_{r},r\in[0,s]\}$,
which is to say that,  
for any positive integer $K$   and
$s_{1},...,s_{K}\leq s$,
 the set of random variables  
$$
\{I_{\sigma_{n}\leq s_{k}}\,
(=I_{\pi_{s_{k}}\geq n}):n=1,2,...,  \;\;  k=1,2,..., K\}
$$
   and 
$\pi_{t}-\pi_{s}$ are independent. (That $\pi_{t}$
introduced in this way possesses the above listed properties
is often put under the rug. For the shortest check, we know,
 see Exercise
2.3.8  and the hint to it  in \cite{Kr02}).

Then take a function $f(t,x,y)$ of   class
$ B_{c}((0,T),C^{\infty}_{0}(\bR^{2}))$
and
  for each $  \omega \in \Omega  $ and $y\in\bR$ solve the equation
 \begin{equation} \label{para}
\partial_{t}u(t,x,y, \omega)=D^{2}_{x}u(t,x,y, 
\omega)+f(t,x,y-h\pi_{t}(\omega))  
\end{equation}
with zero initial data, where $h\in\bR$ is a
parameter. 
 As usual in probability theory in the sequel,  
 more often than not, we
 do not indicate  the dependence on $\omega$.  
 Moreover, we also drop the dependence on $h $ in the sequel. 
 By the above, there exists
   a unique 
   solution
$u(t,x,y)$, depending on $y$ and $\omega$ as parameters,
such that
estimates \eqref{1.23.1}, \eqref{1.22.2},
\eqref{1.22.1}, and
\eqref{1.24.6} hold  
for each $\omega$ and $y\in\bR$
if we replace $u(t,x)$
and $f(t,x)$ with $u(t,x,y)$
and $f(t,x,y-h\pi_{t})$, respectively. Furthermore, 
since
$f\in B_{c}((0,T),C^{\infty}_{0} (\bR^{2}))$,
 $u(t,x,y)$ is uniformly continuous with respect to $y$
uniformly with respect to $\omega,t$,   $h$, and $x$
(cf. the proof of Lemma \ref{lemma 1.25.1}).

By considering $u(t,x,y+h\pi_{t})$ on each interval
$[\sigma_{n},\sigma_{n+1})$ on which $h\pi_{t}$
is constant, one easily derives that
  $u(t,x,y+h\pi_{t})$ satisfies
\begin{equation}
                                      \label{1.23.2a}
 u(t,x,y+h\pi_{t})=
\int_{0}^{t}[D^{2}_{x}u(s,x,y+h\pi_{s})+f(s,x,y)]\,ds
+\int_{(0,t]}g(s,x,y)\,d\pi_{s}
\end{equation}
$$
 = \int_{0}^{t}[D^{2}_{x}u(s,x,y+h\pi_{s})+f(s,x,y)]\,ds
+ \sum_{\sigma_n \le t } g(\sigma_n,x,y),  
$$
where
\begin{equation}
\label{g12}
g(s,x,y)=u(s,x,y+h+h\pi_{s-})-
u(s,x,y+h\pi_{s-}) 
\end{equation}
is the jump of the process $u(t,x,y+h\pi_{t})$
as a function of $t$ at moment $s$ if $\pi_{t}$ has a jump at $s$.

  Here $\pi_{s-} = \lim_{t \uparrow s } \pi_t$,
 $s >0$.  For instance, if $t \in [\sigma_1, \sigma_2)$ we have
\begin{gather*}
u(t,x,y+h\pi_{t})= \int_{0}^{\sigma_1}[D^{2}_{x}u(s,x,y)+f(s,x,y)]\,ds
\\
+ u(\sigma_1,x,y+h)- u(\sigma_1,x,y)
+
\int_{\sigma_1}^{t}[D^{2}_{x}u(s,x,y+h)+f(s,x,y)]\,ds.
\end{gather*}
The next result follows from the   theory of  stochastic 
integrals against     
$\pi_t - \lambda t$ (see Exercise 2.7.8 in \cite{Kr02}).  
We provide a direct and self-contained proof
although a more general situation   will be   encountered
in Lemma \ref{lemma 2.29.2} and treated in a more sophisticated
way.

\begin{lemma}
                                     \label{lemma 1.23.1}
 
For $g$ introduced in \eqref{g12}  and $t\leq T$
we have
$$
 E \int_{(0,t]}g(s,x,y)\,d\pi_{s} 
=\lambda \int_{0}^{t}[v(s,x,y+h)-v(s,x,y)]\,ds,  
$$
where   
\begin{equation} 
                                               \label{v11}
v(t,x,y):=Eu(t,x,y+h\pi_{t}).
\end{equation}
\end{lemma}

Proof. First assume that $t=1$.
 Fix $x$ and $y$ and set
$ g(s)=g(s,x,y)$. 
The function $g$ is bounded on $\Omega\times(0,T)$
and $\pi_{s-}$ is left-continuous with respect to $s$. 
Therefore,
if we define
$$ 
g_{n}(s) = g(k2^{-n})=u(k2^{-n},x,y+h+h\pi_{k2^{-n} - })-
u(k2^{-n},x,y+h\pi_{k2^{-n} -  })
$$
 for $s\in(k2^{-n},(k+1)2^{-n}]$,
$k=0,1,...$, 
then $g_{n}(s)\to  g(s)$ as $n\to\infty$
for any $s\in(0,t]$ and $\omega$,
and
$$
\xi_{n} :=\int_{(0,1]}g_{n}(s)\,d\pi_{s}\to
\int_{(0,1]}g (s)\,d\pi_{s}=:\xi 
$$
 for any  $\omega$. 
By the dominated convergence theorem $E\xi_{n}  
\to E\xi  $.

Next, observe that
\begin{equation}
                                      \label{1.23.3}
E\xi_{n}  =\sum_{k=0}^{2^{ n}-1 } Eg(k2^{-n} )
(\pi_{ (k+1)2^{-n} }-\pi_{ k2^{-n} }).
\end{equation}
Here, owing to the way $g$
was constructed, $g(k2^{-n} )$ is uniquely defined once
we know 
  the values of the random variables
$ I_{\sigma_{i}\leq t  }$ for all
$i=1,2,...$,  and all $t\leq k2^{-n}$,
 and,
 as we have said,
the increments of $\pi_{s}$ after time $k2^{-n}$  
are independent of those 
random variables.
Hence, the expectations of the products 
on the right in \eqref{1.23.3} are equal to the products
of expectations, and since $E\pi_{t}=\lambda t$,
 we conclude, that
$$
E\xi_{n}(t) =\lambda E\sum_{k=0}^{2^{ n}-1 } g(k2^{-n} )
2^{-n}=
 \lambda E\int_{0}^{1}
g_{n}(s)\,ds  
$$
$$
  \to \lambda E\int_{0}^{1}g(s)\,ds=
\lambda \int_{0}^{1}Eg(s)\,ds .
$$
Since, for any $s>0$, we have $\pi_{s}=\pi_{s-}$  (a.s.),
it holds that
$$
Eg(s)=v(s,x,y+h)-v(s,x,y).
$$ 
We have thus proved the lemma if $t=1$.
If it is not, one should just replace above $k2^{-n}$
and $(k+1)2^{-n}$ with $tk2^{-n}$
and $t(k+1)2^{-n}$.
This proves the lemma.  \qed

By taking expectations of both sides of   \eqref{1.23.2a}  
we now obtain the existence part in  
 the following result.

\begin{lemma}
                                         \label{lemma 1.23.2}
Let $f\in B_{c}((0,T),C^{\infty}_{0}(\bR^{2}))$,  
 $h\in \bR$ and $\lambda>0$.
 Then there exists a unique 
continuous function $v(t,x,y)$, $t\in[0,T]$, $x,y\in\bR$,
satisfying the equation
\begin{equation}
                                      \label{1.23.05}
\partial_{t}v(t,x,y)=D^{2}_{x}v(t,x,y)+
\lambda[v(t,x,y+h)-v(t,x,y)]+
f(t,x,y)
\end{equation}
for $t\in(0,T)$, $x,y\in \bR$, with zero initial condition
and such that  
$v(t,\cdot,y)\in C^{2+\alpha}(\bR)$ for any 
$t\in(0,T)$, $ y\in \bR$ and
\begin{equation}
                                                \label{h11}
\sup_{(t,y)\in
[0,T]\times\bR }\|v(t,\cdot,y)\|_{C^{2+\alpha}(\bR)}
\leq N_{0}(T,\alpha) \sup_{
(t,y)\in
(0,T)\times\bR }
\|f(t,\cdot,y)\|_{C^{\alpha}(\bR)}.
\end{equation}
Furthermore,      
$$
\sup_{ (t,z)\in
[0,T]\times\bR^{2} }|v (t,z)|
\leq  T \sup_{(t,z)\in (0,T)\times \bR^{2}}|f (t,z) |,
$$
\begin{equation}
                                                \label{6.30.1}
\sup_{(t,y)\in
[0,T]\times\bR }[D_{x}^{2}v(t,\cdot,y)]_{C^{\alpha}(\bR)}
\leq N_{0}(\alpha) 
\sup_{(t,y)\in
(0,T)\times\bR }[f(t,\cdot,y)]_{C^{\alpha}(\bR)},
\end{equation}
$$
\|D^{2}_{x}v\|^{p}_{L_{p}((0,T)\times\bR^{2})}
\leq N_{p}\|f\|^{p}_{L_{p}((0,T)\times\bR^{2})}
$$
\(where $N_0(T, \alpha),$ $N_0(\alpha)$ and $N_p$ are 
the same as in  \eqref{1.23.1}, \eqref{1.22.1} and 
\eqref{1.24.6}\). 
\end{lemma}

Proof. Uniqueness follows from  
\eqref{1.22.2}  if  $  \lambda T \leq 1/4 $ 
and extends beyond $1/(4\lambda)$ by steps of size 
$1/(4\lambda)$. 

All claimed estimates, apart from the last one, are obtained in the same manner
following the example:
$$
\sup_{y \in \bR}[D_{{x}}^{2}v(t,\cdot,y)]_{C^{\alpha}(\bR)}
\leq \sup_{y\in \bR}E[D_{{x}}^{2}u(t,\cdot,y+h\pi_{t})]_{C^{\alpha}(\bR)},
$$
 where,    for any $t\leq T$ and $\omega$,
$$
\sup_{y\in \bR}[D_{{ x}}^{2}u(t,\cdot,y+h\pi_{t})]_{C^{\alpha}(\bR)} =
\sup_{y\in \bR}[D_{{x}}^{2}u(t,\cdot,y)]_{C^{\alpha}(\bR)}
$$
$$
\le N_{0}(\alpha)\sup_{y \in \bR,s< t}
[f(s,\cdot,y - h \pi_s)]_{C^{\alpha}(\bR)}
\le N_{0}(\alpha) \sup_{y \in \bR,s< T}
[f(s,\cdot,y )]_{C^{\alpha}(\bR)},
$$
which leads to \eqref{6.30.1}.

The last   $L_p$-estimate  is obtained by replacing the above sups
with integrals:
$$
\int_{0}^{T}\int_{\bR^{2}}|D^{2}_{x}v(t,x,y)|^{p}\,dydxdt
\leq E\int_{0}^{T}\int_{\bR^{2}}
|D^{2}_{x}u(t,x,y+h\pi_{t})|^{p}\,dydxdt
$$
$$
=E\int_{0}^{T}
\int_{\bR^{2}}|D^{2}_{x}u(t,x,y )|^{p}\,dydxdt
\leq N_{p}
E\int_{0}^{T}
\int_{\bR^{2}}|f(t,x,y-h\pi_{t} )|^{p}\,dydxdt
$$
$$
=N_{p}
E\int_{0}^{T}\int_{\bR^{2}}
|f(t,x,y )|^{p}\,dydxdt
=N_{p}
 \int_{0}^{T}
\int_{\bR^{2}}|f(t,x,y )|^{p}\,dydxdt.
$$
The lemma is proved.   \qed

We succeeded in adding  
in the right-hand side of \eqref{1.23.5} 
the first-order difference   {\em without changing constants
in our estimates.}  

In our next step, we do with \eqref{1.23.05} almost
the same thing as with \eqref{1.23.5} adding another 
finite difference.
Namely, we introduce $v(t,x,y)$ depending also on 
$\omega$ as a unique solution of
$$
\partial_{t}v(t,x,y)=D^{2}_{x}v(t,x,y)+
\lambda[v(t,x,y+h)-v(t,x,y)]+
f(t,x,y+h\pi_{t})
$$
with zero initial condition.
Then by just repeating the above computations, we see
that
$$
w(t,x,y):=Ev(t,x,y-h\pi_{t})
$$
satisfies
$$
\partial_{t}w(t,x,y)=D^{2}_{x}w(t,x,y)
$$
\begin{equation}
                                      \label{1.24.5}
+
\lambda[w(t,x,y+h)-2w(t,x,y)+w(t,x,y-h)]+
f(t,x,y)
\end{equation}
and admits the same estimates as in Lemma
\ref{lemma 1.23.2}.

Then we take $\lambda=h^{-2}$ in \eqref{1.24.5}
and let $h\downarrow0$. With some extra work,
to be presented later (see the proof of Lemma
\ref{lemma 1.25.1}), one can show that
the solutions   $w= w_{h}$ of \eqref{1.24.5}
with $\lambda=h^{-2}$ converge to a function
$v(t,x,y)$, that is infinitely
differentiable with respect to $(x,y)$
for any $t$ with any derivative
bounded on $[0,T]\times \bR^{2}$,
is continuous in $[0,T]\times \bR^{2}$,
equals zero for $t=0$,
  satisfies
\begin{equation}
                                        \label{6.30.2}
\partial_{t}v(t,x,y)=\Delta v(t,x,y)+f(t,x,y)
\end{equation}
in $(0,T)\times\bR^{2}$ and for which all
the estimates in Lemma \ref{lemma 1.23.2}
hold true with the same constants.  
 
One knows that bounded continuous
in $[0,T]\times \bR^{2}$
solution of \eqref{6.30.2}
having continuous second-order derivatives 
with respect to $(x,y)$ and vanishing at
$t=0$
 are unique, and we conclude
that, for any such solution
 the estimates in Lemma \ref{lemma 1.23.2}
hold true.

Take a unit vector $l_{1}\in\bR^{2}$ and a  unit vector
$l_{2}\in\bR^{2}$ orthogonal to $l_{1}$. Let $S$ 
be an  orthogonal 
transformation of   $\R^{2}$  such that $Se_{i}=l_{i}$,
$i=1,2$, where $e_{1},e_{2}$ is the standard
basis in $\bR^{2}$, and set $f(t,xe_{1}+ ye_{2}) =f(t,x,y)$,
$v(t,xe_{1}+ ye_{2})=v(t,x,y)$,
$$
S(x,y)=xl_{1}+yl_{2},\quad
g(  t ,x,y)=f(t,S(x,y)),\quad w(t,x,y)=v(t,S(x,y)).
$$
 Since the Laplacian is rotation invariant, we have
$$
\partial_{t}w(t,x,y)=\Delta w(t,x,y)+g(t,x,y)
$$
and, since $g$ is as good as $f$, we conclude 
by defining
$$
K=\sup_{(t,y)\in  (0,T)\times\bR}\,\,
\sup_{  x_{1},x_{2}\in\bR, x_1 \not = x_2 }
\frac{|g(t,x_{1},y)-g(t,x_{2},y)|}
{|x_{1}-x_{2}|^{\alpha}}
$$
that
\begin{equation}
                                          \label{8.31.1}
\sup_{(t,y)\in[0,T]\times\bR} \,\, \sup_{
 x_{1} \not = x_{2} }
\frac{|D^{2}_{x}w(t,x_{1},y)-D^{2}_{x}w(t,x_{2},y)|}
{|x_{1}-x_{2}|^{\alpha}}
\leq N_{0}(\alpha)K.
\end{equation}
Observe that, as is easy to see,
 $$
D^{2}_{x}w(t,x,y)=(D^{2}_{l_{1}}v)(t,S(x,y))
=(D^{2}_{l_{1}}v)(t,x l_{1} +y l_{2}) ,
$$ 
where  
 $$
  D_{l}^{2}=l^{i}l^{j}D_{ij} \;\; \text{and}  \;\; 
D_{i}=\partial/\partial x^{i}, \quad D_{ij}=D_{i}D_{j}.
$$
Therefore, the left-hand side  of \eqref{8.31.1}
equals
$$
\sup_{(t,y)\in[0,T]\times\bR}\,\,\sup_{x ,\nu,\mu\in\bR ,
  \mu \not = \nu   }
\frac{|D_{l_{1}}^{2}v(t,\mu l_{1}+xl_{1}+yl_{2})-
D_{l_{1}}^{2}v(t,\nu l_{1}+xl_{1}+yl_{2})|}
{|\mu-\nu|^{\alpha}} 
$$
$$
=\sup_{(t,z)\in[0,T]\times\bR^{2}}\,\,\sup_{     \mu \not = \nu  }
\frac{|D_{l_{1}}^{2}v(t,\mu l_{1}+z)-
D_{l_{1}}^{2}v(t,\nu l_{1}+z)|}
{|\mu-\nu|^{\alpha}}.
$$
Similarly the right-hand side of \eqref{8.31.1}
is transformed and we get that 
for any (actually, only one) bounded continuous
in $[0,T]\times \bR^{2}$
solution $v$ of \eqref{6.30.2}
having continuous second-order derivatives 
with respect to $(x,y)$ and vanishing at
$t=0$ and any unit vector $l\in\bR^{2}$
$$
\sup_{(t,z)\in[0,T]\times\bR^{2}}\,\,\sup_{   \mu \not = \nu  }
\frac{|D_{l}^{2}v(t,\mu l +z)-D_{l}^{2}v(t,\nu l +z)|}
{|\mu-\nu|^{\alpha}}
$$
\begin{equation}
                                                       \label{10.9.2}
\leq N_{0}(\alpha)\sup_{(t,z)\in(0,T)
\times\bR^{2}}\,\,\sup_{  \mu \not = \nu  }
\frac{|f(t,\mu l +z)-f(t,\nu l +z)|}
{|\mu-\nu|^{\alpha}}.
\end{equation}
Also, since the Jacobian
of the above $S(x,y)$ equals one,
 for any unit vector $l\in\bR^{2}$
\begin {equation}                            \label{sti1}
\int_{0}^{T}\int_{\bR^{2}}|D^{2}_{l}v(t,z)|^{p}\,dzdt
\leq N_{p}
\int_{0}^{T}\int_{\bR^{2}}|f(t,z)|^{p}\,dzdt.
\end{equation}

\mysection{Multidimensional second-order parabolic
equations}
                                      \label{section 6.29.2}

\begin{theorem}
                                \label{theorem 1.25.1}
Let $a(t)=(a^{ij}(t))$ be a $d\times d$ symmetric
matrix-valued   Borel measurable function on $(0,T)$
such that 
\begin{equation}
                                        \label{1.25.12}
a^{ij}(t)\lambda^{i}\lambda^{j}\geq|\lambda|^{2}
\end{equation}
for all $t\in(0,T)$ and $\lambda\in\bR^{d}$  
and
\begin{equation}
                                       \label{9.29.1}
\int_{0}^{T}\tr a(t)\,dt<\infty.
\end{equation}   
Then for any $f\in B_{c}((0,T),C^{\infty}_{0} (\bR^{d}))$
there exists a unique    continuous in $[0,T]\times\bR^{d}$
solution $u(t,x)$ of the equation
\begin{equation}
                                        \label{1.25.1}
\partial_{t}u(t,x)=a^{ij}(t)D_{ij}u(t,x)+f(t,x)
\end{equation}
in $(0,T)\times\bR^{d}$ with zero initial data
such that, 
 for any $t\in[0,T]$, $u(t,\cdot)\in C^{2+\alpha}(\bR^{d})$
and,
for any $i.j=1,...,d$ and unit vector $l\in\bR^{d}$,
we have  
\begin{equation}
                                        \label{1.25.2}
 \sup_{ (t,x) \in 
[0,T]\times\bR^{d} }|u (t,x) |
\leq  T \sup_{ (t,x) \in  (0,T)\times \bR^{d}}|f (t,x)|,
\end{equation}
\begin{equation}
                                        \label{1.25.3}
\sup_{ t \in 
[0,T] }[D_{ij} u(t,\cdot )]_{C^{\alpha}(\bR^{d})}
\leq N'(\alpha)N_{0}(\alpha) 
\sup_{ t \in 
(0,T) }[ f(t,\cdot )]_{C^{\alpha}(\bR^{d})},
\end{equation}
\begin{equation}
                                        \label{1.25.5}
\sup_{  (t,x) \in 
[0,T]\times \bR^{d}}[D_{l}^{2} u(t,x+l\,\cdot )]_{C^{\alpha}(\bR )}
\leq  N_{0}(\alpha) 
\sup_{  (t,x) \in 
(0,T)\times \bR^{d} }[  f(t,x+l\,\cdot)]_{C^{\alpha}(\bR )},
\end{equation}
\begin{equation}
                                        \label{1.25.4}
\|D^{2}_{l}u\|^{p}_{L_{p}((0,T)\times\bR^{d})}
\leq N_{p}\|f\|^{p}_{L_{p}((0,T)\times\bR^{d})},
\end{equation}
where   $N_{0}(\alpha)$,
$N_{p}$ are the constants from Section \ref{section 1.25.1} 
 \(see \eqref{1.22.1} and \eqref{1.24.6}\) 
and $N'(\alpha)$ is a constant specified
in Lemma \ref{lemma 1.25.2}.
\end{theorem}

We see, in particular, that  the $L^{1}$-norms
of $a^{ij}(t)$ do  not influence the constants in the estimates.

\begin{lemma}
                                     \label{lemma 1.25.1}
The assertions of Theorem \ref{theorem 1.25.1},
apart from \eqref{1.25.3}, hold true if $a^{ij}=\delta^{ij}$.
\end{lemma}

Proof. We proceed by induction on $d$.
Assume that the lemma is true for a particular $d$
and repeat the construction in     Lemma \ref{lemma 1.23.2}  treating $x$ there as a point in
$\bR^{d}$ and replacing $D^{2}_{x}$ with
the Laplacian $\Delta_{x}$ 
in $\bR^{d}$.
  Then, under the assumption that we are given
$f(t,x,y)$, $t\in(0,T),x\in\bR^{d},y\in\bR$,
which is of class $B_{c}((0,T),C^{\infty}_{0}(\bR^{d+1}))$,
 we arrive at the conclusion that, 
for any $h>0$,  the equation
$$
\partial_{t}u_{h}(t,x,y)=\Delta_{x}u_{h}(t,x,y)
+
f(t,x,y)
$$
\begin{equation}
                                      \label{1.25.6}
+
h^{-2}[u_{h}(t,x,y+h)-2u_{h}(t,x,y)+u_{h}(t,x,y-h)], 
\end{equation}
where $t\in(0,T),x\in\bR^{d},y\in\bR$, with zero initial
condition has a unique continuous in $[0,T]\times\bR^{d+1}$
solution $u_{h}(t,x,y)=u_{h}(t,z)$, where $z=(x,y)$,
 such that   
\begin{equation}
                                        \label{1.25.7}
 \sup_{ (t,z) \in
[0,T]\times\bR^{d+1} }|u_{h} (t,z) |
\leq  T \sup_{ (t,z) \in (0,T)\times \bR^{d+1}}|f(t,z) |,
\end{equation}
$$
\sup_{  (t,z) \in
[0,T]\times \bR^{d+1}}[D_{l}^{2} 
u_{h}(t,x+l\,\cdot,y )]_{C^{\alpha}(\bR )}
$$
\begin{equation}
                                        \label{1.25.8}
\leq  N_{0}(\alpha) 
\sup_{  (t,z) \in
(0,T)\times \bR^{d+1} }
[  f(t,x+l\,\cdot,y)]_{C^{\alpha}(\bR )},
\end{equation}
\begin{equation}
                                        \label{1.25.9}  
\|D^{2}_{l}u_{h}\|^{p}_{L_{p}((0,T)\times\bR^{d+1})}
\leq N_{p}\|f\|^{p}_{L_{p}((0,T)\times\bR^{d+1})},
\end{equation}
where $l$ is any unit vector in $\bR^{d}$.

One can apply the finite-difference operators
with respect to $(x,y)$
of any order to \eqref{1.25.6};   these operators are obtained by 
compositions of the first order difference operators like 
$$
\delta_{r,i}v(z)=r^{-1}[v(z+re_{i})-v(z)],\quad i=1,...,d+1,
$$
where $e_{i}$ is the $i$th basis vector and $r>0$.

Then,
owing to \eqref{1.25.7} and the fact that
any derivative of any order of $f$
is in $B_{c}((0,T),C^{\infty}_{0}(\bR^{d+1}))$, 
we conclude that any finite-difference of any order
of $u_{h}$ is bounded on $\bR^{d+1}$
uniformly with respect to $t$ and $h$.
It follows that $u_{h}$ is infinitely differentiable
with respect to $(x,y)$ and any derivative of any order
is bounded on $[0,T]\times\bR^{d+1}$. Then equation
\eqref{1.25.6} itself (always considered
in the integral form as \eqref{ci1}) shows that 
 these derivatives are Lipschitz continuous
in $t$.
Thus, the family $u_{h}$ is equi-Lipschitz in each compact set of $[0,T] \times \bR^{d+1}$ and the same holds for any derivative with respect to $(x,y)$
of $u_{h}$. 
 
Now by the Arzel\`a-Ascoli theorem
there is a sequence $u_{h_{n}}$, $h_{n}\downarrow 0$,
which converges uniformly on any set
$[0,T]\times\{|(x,y)|\leq R\}$, $R\in(0,\infty)$,
along with any derivative with respect to $(x,y)$
of $u_{h_{n}}$ and $\partial_{t}u_{h_{n}}$.

   Writing \eqref{1.25.6} in the integral form as \eqref{ci1} 
 and passing to the limit   as $n \to \infty$, 
 we conclude that
there exists a continuous function $u(t,x,y)$ 
in $[0,T]\times\bR^{d+1}$, which is infinitely
differentiable with respect to $(x,y)$
with any derivative bounded on $[0,T]\times\bR^{d+1}$;   moreover, the equation
$$
\partial_{t}u(t,x,y)=\Delta_{x,y}u(t,x,y)+f(t,x,y)
$$
holds in $(0,T)\times\bR^{d+1}$ and estimates
\eqref{1.25.7}, \eqref{1.25.8}, and
\eqref{1.25.9}  are valid with $u$ in place of $u_{h}$.

Uniqueness of such solutions is a simple consequence
of   
the maximum principle.
The invariance
of the Laplacian in $\bR^{d+1}$ under rotations
shows that estimates \eqref{1.25.2},
\eqref{1.25.5}, and \eqref{1.25.4} are true
with $\bR^{d+1}$ in place of $\bR^{d}$
for any unit vector $l\in\bR^{d+1}$ 
(cf.  \eqref{10.9.2} and \eqref{sti1}).
The lemma is proved.              \qed

The following lemma shows that \eqref{1.25.3}
follows from \eqref{1.25.5}.

\begin{lemma}
                                    \label{lemma 1.25.2}
Let $u\in C^{2+\alpha}(\bR^{d})$ be   such that, for any unit vector
$l\in\bR^{d}$, we have   
$$
\sup_{ x \in
 \bR^{d}}[D_{l}^{2} u( x+l\,\cdot )]_{C^{\alpha}(\bR )}
\leq 1.
$$
Then there exists a constant $N'(\alpha)$
such that for any $i,j=1,...,d$ we have
$$
M:=[D_{ij}u]_{C^{\alpha}(\bR^{d})}\leq N'(\alpha).
$$
\end{lemma}

Proof. We use the method of proof    
which we learned
from M. Safonov.
 Let $T_{x_{0}}(x)$ denote the second-order Taylor
polynomial of $u$ centered at $x_{0}$. Then by the 
mean-value theorem for any unit vector $l\in\bR^{d}$ and $t\geq0$
$$
|u(x_{0}+tl)-T_{x_{0}}(x_{0}+tl)| =   (1/2)t^{2}
|D^{2}_{l}u(x_{0}+\theta l) -D^{2}_{l}u(x_{0})|
\leq (1/2)t^{2+\alpha},
$$
where $\theta\in(0,t)$. It follows that
for any $r\in(0,\infty)$ and $x_{0}\in\bR^{d}$
there exists a quadratic polynomial $p(x)$ such that
$$  
|u(x)-p(x)|\leq (1/2)r^{2+\alpha}  
$$
in $B_{r}(x_{0})=\{x:|x-x_{0}|<r\}$.

Observe that by the mean-value theorem, for $h>0$,
$$  
|D_{ij}u(x)-\delta_{h,i}\delta_{h,j}u(x)|\leq M(2h)^{\alpha}.
$$
Next, take $x_{1},x_{2}\in \bR^{d}$,
choose $h=\varepsilon|x_{1}-x_{2}|$, where
$\varepsilon$ is such that
$$
(2\varepsilon)^{\alpha}=1/4,
$$
and observe, that
if $r=
|x_{1}-x_{2}|+2h$, then all six points
$x_{k},x_{k}+he_{i},x_{k}+he_{i}+he_{j}$, $k=1,2$,
can be encompassed by a ball of radius $r$ 
(centered at $x_1$). 
By the above, for an appropriate quadratic
polynomial $p$ 
  (we use the fact that $\delta_{h,i}\delta_{h,j}p$ is constant 
since it is a 
  second-order difference of a quadratic polynomial)  
$$
|D_{ij}u(x_{1})-D_{ij}u(x_{2})|\leq
(1/2)M|x_{1}-x_{2}|^{\alpha}
$$
$$
+|\delta_{h,i}\delta_{h,j}(u-p)(x_{1})
-\delta_{h,i}\delta_{h,j}(u-p)(x_{2})|,
$$
where the last term is less than
$$
|\delta_{h,i}\delta_{h,j}(u-p)(x_{1})|
+|\delta_{h,i}\delta_{h,j}(u-p)(x_{2})|
$$
$$
\leq  3r^{2+\alpha}h^{-2}
= 
 3(1+2\varepsilon)^{2+\alpha}
\varepsilon^{-2}|x_{1}-x_{2}|^{\alpha}.
$$
The arbitrariness of $x_{1}$ and $x_{2}$ now yields
the desired result with
$$
N' (\alpha )=6(1+2\varepsilon)^{2+\alpha}\varepsilon^{-2}.
$$
The lemma is proved.     \qed

 In the sequel, given a unit vector $l \in \bR^d$, we 
denote by $l l^*$ the $d\times d$ matrix with entries
$l^{i}l^{j}$.

\begin{lemma}
                                       \label{lemma 1.25.4}
Let the assertions of Theorem
\ref{theorem 1.25.1} be true for a given
$a(t)$ satisfying the assumptions of the theorem
and such that it is continuous. Let
$\nu(t)$ be a real-valued  continuous  function
on $[0,T]$ and $l\in\bR^{d}$ be a unit vector.
Then the assertions of Theorem
\ref{theorem 1.25.1} hold true for
$a(t)+\nu^{2}(t)ll^{*}$, as well, with the same
constants in the estimates
(hence the constants are independent of $\nu(t)$ and $l$).
\end{lemma}

Proof. Introduce  
$$
b_{t}=l\int_{(0,t]}\nu(s)\,d\pi_{s}\quad
\big(= l \sum_{\sigma_{n}\leq t} \nu (\sigma_{n})=
l \sum_{s\leq t} \nu (s)(\pi_{s}-\pi_{s-})\big).
$$
Observe that for $0\leq s\leq t<\infty$ 
\begin{equation}
                                               \label{1.24.2}
 E(b_{t}-b_{s})=  \lambda   l \int_{s}^{t}  \nu(r)\,dr 
\end{equation}
(which is easily proved if $\nu$ is piece-wise constant,
and then extended to continuous   $\nu$
by standard arguments, 
cf. the proof of Lemma \ref{lemma 1.23.1}). 

 Then take a function $f(t,x)$ of class
$ B_{c}((0,T),C^{\infty}_{0}(\bR^{d}))$
and
  for each $\omega$  solve the equation
$$
\partial_{t}u(t,x )=a^{ij}(t)D_{ij}u(t,x )+f(t,x -hb_{t})
$$
with zero initial data, where $h\in\bR$ is a parameter.
   In the sequel we drop the dependence on $h$. 
 By assumption, there exists
   a unique continuous in   $[0,T]\times\bR^d$  solution
$u(t,x )$ depending on  $\omega$ as parameter 
such that
estimates \eqref{1.25.2}, \eqref{1.25.3},
\eqref{1.25.5}, and
\eqref{1.25.4} hold  
for each $\omega$  
if we replace 
  $f(t,x)$ with  
   $f(t,x - h b_{t})$  (which, by the way, does
not affect the right-hand sides of  these estimates).  Furthermore, 
since
$f\in B_{c}((0,T),C^{\infty}_{0} (\bR^{d}))$,  
 $u(t,x )$ is uniformly continuous with respect to $x$
uniformly with respect to $\omega,t$,   and $h$ (cf. the proof
of Lemma \ref{lemma 1.25.1}).

By considering $u(t,x )$ on each interval
$[\sigma_{n},\sigma_{n+1})$ on which $ \pi_{t}$,
and hence $b_{t}$,
are constant, one easily derives that
  $u(t,x+hb_{t} )$ satisfies
\begin{equation}
                                      \label{1.23.2}
u(t,x +hb_{t})=
\int_{0}^{t}[a^{ij}(s)D_{ij}u
 (s,x +hb_{s})+f(s,x )]\,ds
+\int_{(0,t]} g(s,x)\,d\pi_{s},
\end{equation}
where
$$
g(s,x):=u(s,x +hl\nu(s)+hb_{s-})-  
u(s,x +hb_{s-}).
$$  
By introducing
$$
g_{n}(s,x)=u(k2^{-n},x +hl\nu(k2^{-n})+hb_{k2^{-n} -})-  
u(k2^{-n},x +hb_{k2^{-n}-})
$$
  for $s\in(k2^{-n},(k+1)2^{-n}]$,
$k=0,1,...$,  using the continuity
of $\nu(t)$ and \eqref{1.24.2},
and repeating the proof of
Lemma \ref{lemma 1.23.1}, we arrive at the conclusion
that
$$
E\int_{(0,t]} g(s,x)\,d\pi_{s}=\lambda
\int_{0}^{t} [v(s,x+hl\nu(s))-v(s,x))]\,ds,
$$
where
$$
v(t,x)=Eu(t,x+hb_{t}).
$$
 
Then \eqref{1.23.2} yields
$$
\partial_{t}v(t,x)=a^{ij}(t)D_{ij} v
 (t,x )+\lambda [v(t,x+hl\nu(t))-v(t,x)]+f(t,x ).
$$

As in Section \ref{section 1.25.1}, $v$ is a unique
  solution
of this equation for which all estimates
claimed in the theorem hold true.

After that we solve
$$
\partial_{t}w(t,x)=   a^{ij}(t)D_{ij} w  
 (t,x )+\lambda [w(t,x+hl\nu(t))-w(t,x)]+f(t,x+hb_{t} )
$$
and repeat the end of Section \ref{section 1.25.1}
to conclude that for each $h>0$ there exists
a continuous function $u_{h}(t,x)$ on $[0,T]
\times\bR^{d}$, which is a unique solution of
$$
\partial_{t}u_{h}(t,x)=   a^{ij}(t)D_{ij} u_h 
 (t,x )+f(t,x)
$$
$$
+h^{-2}[u_{h}(t,x+hl\nu(t))-2u_{h}(t,x)
+u_{h}(t,x-hl\nu(t))]
$$
in $(0,T)\times\bR^{d}$ with zero initial condition
and for which all estimates
claimed in the theorem hold true. 

As in the proof of Lemma \ref{lemma 1.25.1},
a subsequence $u_{h_{n}}$ converges to
the function we are after. The lemma is proved.   \qed

{\bf Proof of Theorem \ref{theorem 1.25.1}}. 
  Uniqueness    is easily derived from 
  the maximum principle. 
(Just in case,
if the reader sees any obstacle in the fact that
$a^{ij}$ may be unbounded, have in mind that a trivial
time change  (i.e.,  $u(t,x) = v(\int_0^t \tr a(s)ds, x)$)   reduces the general situation to the
one with $\tr a(t)\equiv1$. Actually, after the time change
the new matrix may degenerate, but this is not
an obstacle for the maximum principle for
parabolic equations to hold, see, for instance,
Theorem 4.1 of \cite{KP}. Also see Corollary 3.6 there.)
 To prove the existence of solutions,
  by having in mind a simple passage to the limit
(we say more about this in  Theorem
\ref{theorem 9.17.1}
 and its proof in Section \ref{section 6.29.6}
or send the reader to the end of the present proof)
and approximating $a(t)$ by $a_{n}(t)=a(t)I_{\tr a(t)\leq n}
+(\delta^{ij})I_{\tr a(t)> n}$,
we may assume that $a(t)$ is bounded. By the same token
we may assume that there exists a constant
$\varepsilon>0$ such that  
\begin{equation}
                                              \label{1.25.11}
a^{ij}(t)\lambda^{i}\lambda^{j}\geq  
(1+2\varepsilon) |\lambda|^{2} 
\end{equation}
for all $t\in(0,T)$ and $\lambda\in\bR^{d}$.

Then
for the matrix $\hat a(t)=(\hat a^{ij}(t))
=(a^{ij}(t)-\delta^{ij})$ we have
 $$
\hat a^{ij}(t)\lambda^{i}\lambda^{j}\geq  2
 \varepsilon  |\lambda|^{2},
$$
for all $t\in(0,T)$ and $\lambda\in\bR^{d}$. 
By assumption $\text{tr\,} \hat a(t)$ is also bounded,
so that $a(t)$ takes values in a closed subset $\Gamma$ of
the set $S( M)$ of symmetric $d\times d$-matrices
$a$ such that
$$
a^{ij}\lambda^{i}\lambda^{j} > \varepsilon 
 |\lambda|^{2},\quad
 \lambda\ne0, \quad
\text{tr\,} a<M.  
$$

One knows that
 there exist
  $n\in\{1,2,...\}$, vectors $l_{1},...,l_{n}\in\bR^{d}$,
and real-analytic  real-valued
functions  $\nu_{1}(a),...,\nu_{n}(a)$ on $S( M)$,
such that for $a\in\Gamma$ it holds that
$$
a=\sum_{k=1}^{n}\nu^{2}_{k}(a)l_{k}l_{k}^{*}
$$
(for instance, see Section 1 in \cite{Kr08}).
In particular,
\begin{equation}
                                              \label{1.25.10}
 a(t)=(\delta^{ij})+\sum_{k=1}^{n}\nu^{2}_{k}(t)l_{k}l_{k}^{*},
\end{equation}
where $\nu_{k}(t)=\nu_{k}(\hat a(t))$.
The functions $\nu_{k}(t)$ are continuous
if $a(t)$ is continuous, and, therefore,
by using Lemma \ref{lemma 1.25.1}
and an obvious induction on the number of terms
in \eqref{1.25.10} along with Lemma \ref{lemma 1.25.4}
we conclude that the theorem holds true
under the additional assumptions that $a(t)$
is continuous and \eqref{1.25.11} holds.

To abandon the continuity assumption,
we find uniformly bounded smooth $a_{n}(t)$, $n=1,2,...$,
 satisfying
\eqref{1.25.12}, such that $a_{n}(t)\to a(t)$
as $n\to\infty$
for almost all $t$.

  We extend $a$ to the whole $\bR$ by setting $a(t) = a(T/2)$, if $t \ge T$ or  $t \le 0.$ Then we consider standard mollifiers $(\rho_n) \subset C_0^{\infty}(\bR)$ and introduce the matrices $a_n(t) = 
(a_n^{ij}(t))$, 
$$
a_n^{ij} (t) = (a^{ij} * \rho_n)(t),\;\; t \in \bR.
$$
It is clear that each $a_n(t)$ is symmetric and non-negative and depends continuously on $t$; 
 moreover  
$$
\sup_{t \in \bR} \tr a_n(t) 
\le \sup_{t \in (0,T)} \tr a(t)
$$ and 
$$
 a^{ij}_n(t)\lambda^{i}\lambda^{j} \ge  |\lambda|^{2},\;\;\; t \in \bR,
\lambda\in\bR^{d}.
$$
Let us consider solutions $u_n$ of 
\begin{equation}
\label{dee}
u_n (t,x) = \int_0^t a^{ij}_n(s) D_{ij} u_n (s,x) ds + 
\int_0^t f (s,x) ds,
\end{equation} 
the ones obtained
according to the first part of the proof.

We can use  estimates \eqref{1.25.2}, \eqref{1.25.3}, and \eqref{1.25.5} 
with $u$ replaced by $u_n$.   Moreover, using also \eqref{dee} we deduce that 
the family $u_{n}$ is equi-Lipschitz in each compact set of 
$[0,T] \times   \bR^d $; the same holds for any derivative with respect to $x$
of $u_{n}$.
 
By the Arzel\`a-Ascoli theorem
there is a subsequence which we still denote by $u_{n}$
which converges uniformly on any set
$[0,T]\times\{|x|\leq R\}$, $R\in(0,\infty)$,
along with any derivative with respect to $x$
of $u_{n}$.

Passing to the limit  as $n \to \infty$ in \eqref{dee} 
 we conclude that
there exists a continuous function $u(t,x)$ 
in $[0,T]\times\bR^{d}$, which is infinitely
differentiable with respect to $x$
with any derivative bounded on $[0,T]\times\bR^{d}$. 
 Such function $u$ is a solution to \eqref{1.25.1}. 
 Moreover  estimates
\eqref{1.25.3}, \eqref{1.25.5} and
 \eqref{1.25.4} hold  for~$u$.  \qed

\begin{corollary}
                                    \label{corollary 1.26.1}
Let $u\in C^{\infty}_{0}(\bR^{d})$
and assume that  
 $a(t)$ in Theorem \ref{theorem 1.25.1}  is 
independent of $t$, i.e., $a(t) \equiv a$.  
Set
$$
f=a^{ij}D_{ij}u.
$$
Then for all $i,j=1,...,d$ and unit vector $l\in\bR^{d}$
we have
$$
[D_{ij}u]_{C^{\alpha}(\bR^{d})}\leq N'(\alpha)N_{0}(\alpha)
[f]_{C^{\alpha}(\bR^{d})},
$$
$$
\|D^{2}_{l}u\|^{p}_{L_{p}(\bR^{d})}\leq N_{p}
\|f\|^{p}_{L_{p}(\bR^{d})}.
$$
\end{corollary}

Proof.   Let $T>0$.    The function   $u(t,x):=u(x)t/T
$  is  a unique  bounded 
solution of
$$
\partial_{t}u(t,x)=a^{ij}D_{ij}u(t,x) + 
 g_{T}(t,x)
$$
 with zero initial condition,  
where $g_{T}(t,x)=u(x)/T -f(x)t/T$. By Theorem 
\ref{theorem 1.25.1}
$$
[D_{ij}u]_{C^{\alpha}(\bR^{d})}\leq N'(\alpha)N_{0}(\alpha)
\left([f]_{C^{\alpha}(\bR^{d})}+(1/T)
 [u]_{C^{\alpha}(\bR^{d})}\right),
$$
$$
\int_{\bR^{d}}|D^{2}_{l}u(x) |^{p}\,dx
\int_{0}^{T}(t/T)^{p}\,dt
\leq  N_{p}
\int_{0}^{T}\int_{\bR^{d}} | u(x)/T -(t/T)f(x)|^{p}\,dxdt,
$$
$$
\int_{\bR^{d}}|D^{2}_{l}u(x) |^{p}\,dx
\leq (p+1)N_{p} /T
 \int_{0}^{T}\int_{\bR^{d}}[u(x)/T+(t/T)f(x)]^{p}\,dxdt 
$$

$$
=(p+1)N_{p}\int_{0}^{1}\int_{\bR^{d}}|u(x)/T + s f(x)|^{p}\,dxds,
$$
and our assertions follow after letting $T\to\infty$. \qed

\begin{remark}
                                    \label{remark 1.26.1}
For fixed $T\in(0,\infty)$ denote by $N_{p}(d)$ the least constant $N$
such that 
$$
\|D^{2}_{l}u\|_{L_{p}((0,T) \times \bR^{d})}^{p}
\leq N\|f\|_{L_{p}((0,T) \times \bR^{d})}^{p}
$$ 
for any unit vector $l\in\bR^{d}$,
 $f\in B_{c}((0,T), C^{\infty}_{0}(\bR^{d}))$, and  any
bounded continuous in $[0,T]\times\bR^{d}$
solution $u$ of the equation
\begin{equation}
                                            \label{1.26.1}
\partial_{t}u=\Delta u+f
\end{equation}
in $(0,T)\times\bR^{d}$ with zero initial condition.
 It turns out that
$$
N_{p}(d)=N_{p}(1). 
$$
Indeed, by Theorem  \ref{theorem 1.25.1},
$N_{p}(d)\leq N_{p}(1)$. On the other hand,
let 
$$
\WO^{1,2}_{p}([0,T]\times \bR^{d})
=\{u\in W^{1,2}_{p}([0,T]\times \bR^{d}))
\cap C([0,T],L_{p}( \bR^{d})):u(0,\cdot)=0\}  
$$
 (for the definition of $W^{1,2}_{p}([0,T]\times \bR^{d})$, 
see for instance, page 153 in \cite{K3}). 
We know that the operator
$\partial_{t}-\Delta$ maps $\WO^{1,2}_{p}([0,T]\times \bR^{d})$
onto $L_{p}((0,T)\times\bR^{d}))$ in a one-to-one way
and has a bounded inverse. Furthermore,
  the set $B_{c}((0,T), C^{\infty}_{0}(\bR^{d}))$
is dense in $L_{p}((0,T)\times \bR^{d})$.
It follows that $N_{p}(d)$
is the least constant $N$ such that
for any $u\in \WO^{1,2}_{p}([0,T]\times \bR^{d})$
and unit vector $l\in\bR^{d}$ we have  
$$
\|D^{2}_{l}u\|_{L_{p}((0,T) \times \bR^{d})}^{p}
\leq N
  \|\partial_{t}u-\Delta u\|_{L_{p}((0,T) 
\times \bR^{d})}^{p}.
$$

Now,
let $ u(t,x )$ be a function
of class $\WO^{1,2}_{p}([0,T]\times \bR )$
and let $\zeta(x')=\zeta(x^{2},...,x^{d})$
be a nonzero function of class $C^{\infty}_{0}(\bR^{d-1})$.
Introduce $u_{n}(t,x)=u(t,x^{1}) \zeta(x'/n)$.
By definition, ($D_{11}=\partial^{2}/(\partial x^{1})^{2}$)
$$
\int_{0}^{T}\int_{\bR}|D_{11}u(t,x )|^{p}\,dxdt
\int_{\bR^{d-1}}\zeta^{p}(y/n)\,dy
$$ 

$$
\leq N_{p}(d)\int_{0}^{T}\int_{\bR^{d}}
\big|\zeta( x'/n)[\partial_{t}u(t,x^{1})
-D_{11}u(t,x^{1})]
$$

$$
- n^{-2}
u(t,x^{1})(\Delta \zeta) (x'/n)
\big|^{p}\,dxdt,
$$

$$
\int_{0}^{T}\int_{\bR}|D_{11}u(t,x )|^{p}\,dxdt
\int_{\bR^{d-1}}\zeta^{p}   (y)  \,dy
$$ 
$$
\leq N_{p}(d)\int_{0}^{T}\int_{\bR^{d}}
\big|\zeta ( x'  )[\partial_{t}u(t,x^{1})
-D_{11}u(t,x^{1})]- n^{-2} 
u(t,x^{1})(\Delta \zeta)(x' )
\big|^{p}\,dxdt.
$$
By letting $n\to\infty$ we get
$$
\int_{0}^{T}\int_{\bR}|D^{2}u(t,x )|^{p}\,dxdt
\leq N_{p}(d) \int_{0}^{T}\int_{\bR}|
\partial_{t}u(t,x )
-D^{2}u(t,x)|^{p}\,dxdt
$$ 
and, since this is true for any element   $u$  of
$\WO^{1,2}_{p}([0,T]\times \bR )$, we have
$N_{p}(d)\geq N_{p}(1)$.

\end{remark}

\mysection{General setting. Main results}
                                         \label{section 6.29.3}

Let $W$ be a 
set
consisting of real-valued 
 (Borel) measurable functions 
 $u=u_{t}=u_{t}(x)$ 
on $[0,T]\times\bR^{d}$.
 In Sections \ref{section 1.25.1} and \ref{section 6.29.2}
we only considered bounded solutions. Therefore,
we assume that the elements of $W$ are bounded
and even uniformly bounded as required in Assumption 
\ref{assumption 2.26.3} (i) below.

Let $\cG$ be a commutative group 
of affine volume-preserving transformations 
of $\bR^{d}$.  If $g,h\in\cG$ by $gh$  we mean
the composition of the two transformations.
\begin{remark}
                                          \label{remark 8.24.1}
We draw the reader's attention to the fact that,
since each $g\in\cG$ is   measure-preserving,  its
Jacobian equals one.
\end{remark}
 As usual, if $f(x)$
is a function on $\bR^{d}$ and $g\in\cG$, we define
$(gf)(x)=f(gx)$, where $gx$ is  the image of $x$ 
under mapping $g$.

By $B((0,T),\cG)$ we denote the set of
bounded measurable $\cG$-valued functions on $(0,T)$.

Fix a constant $K\in[0,\infty)$.

\begin{assumption}
                                \label{assumption 2.26.3}
(i) 
For any $u\in W$   we have
  $$ 
\sup_{(t,x ) \in [0,T] \times \bR^d} | u_{t}
(x) | 
  \leq K.
$$

 (ii) (Convexity of $W$.) If $(\Omega,\cF,P)$ is a probability space and $u(\omega)
=u_{t}(\omega,x)$ is an $\cF\times \cB([0,T]\times\bR^{d})$-measurable function such that $u(\omega)\in W$ for any $\omega$,
then  
the function $E [u_t (x)]$ belongs to $ W$   (by $E$ we indicate the expectation with respect to $P$ and by ${\mathcal B}  ([0,T]\times\bR^{d})$ we mean the Borel $\sigma$-field on
$[0,T]\times\bR^{d}$).

(iii) (``Shift'' invariance of $W$.)
For $u \in W$ and any bounded  measurable $\cG$-valued function
$g_{t}$  given on $[0,T]$, the function
 $u_{t} ( g_{t}x)$ is  in $ W$.
 
\end{assumption}

Let $L:=\{L_{t}, t\in(0,T)\}$, be a family of linear  operators
$$
L_{t}:C^{\infty}_{0}(\bR^{d})\to B(\bR^{d})
$$   
 ($B(\bR^{d})$ denotes the  space of
real-valued bounded and Borel functions defined on 
$\bR^d$)  
and take and fix
\begin{equation}
                                               \label{5.8.4}
f\in  B ((0,T)\times\bR^{d} ) ,\quad
u_{0}\in  
   B(\bR^{d}),  
  \end{equation}
 where $B ((0,T)\times\bR^{d} )$ is the set 
of Borel bounded functions
  on $(0,T)\times\bR^{d}$.

\begin{assumption}
                                   \label{assumption 2.27.1}
The couple $(L,f)$ is $W$-regular in the 
following sense.

(i) ($\cG$ and $L$ commute.)
For any $t\in(0,T)$ and $g\in\cG$,   we have $gL_{t}=L_{t}g$.

(ii)  
For any $\zeta\in  C^{\infty}_{0}(\bR^{d})$,
$L_{t}\zeta (x):=(L_{t}\zeta)(x)$ is measurable with respect to 
$(t,x)$ and
$$
 \int_{ [0,T]\times\bR^{d} } |L_{t}\zeta(x)| dt dx <\infty .
$$
 
(iii) 
There is a mapping $B((0,T),\cG)\to W$ mapping every bounded
measurable
$\cG$-valued functions  $h=h_{t},t\in(0,T)$, 
 into $u[h]\in W$
 such that $u=u[h]$ satisfies the
equation 
\begin{equation}
                                               \label{2.28.1}
u_{t}(x)=u_{0}(x)+\int_{0}^{t}[ L^{*}_{r} u_{r}(x) 
+(h_{r}f_{r})(x)]\,dr,\quad t\in[0,T],\; x \in \bR^d,
\end{equation}
in the sense specified below (see \eqref{weak11}). 

(iv) For any  $h',h''\in B((0,T),\cG)$   
and $(t,x)\in[0,T]\times \bR^{d}$,
we have
\begin{equation}                                 \label{fer}
|u_{t}[ h'](x)-u_{t}[  h''](x)|\leq K\int_{0}^{t}
\sup_{y \in \bR^{d}}|f_{r}( h'_{r}y)-f_{r}( h'' _{r}y)| \,dr.
\end{equation}
  \end{assumption}

\begin{remark}
                                         \label{remark 5.5.1}
Assumption \ref{assumption 2.27.1} (iv)
 implies that,  for any $h \in B((0,T),\cG)$, $x\in \bR^{d}$,
and $t\leq s \leq T$ we have
$$
u_{t}[ h  ](x)=u_{t}[ h_{\cdot\wedge s}  ](x).
$$ 
Indeed, it is enough to use 
\eqref{fer} with $h'=h$ and $h'' = h_{\cdot\wedge s}$.  
\end{remark} 

We say that $u\in
W$ satisfies \eqref{2.28.1} if,
for any $\zeta\in C^{\infty}_{0}(\bR^{d})$ and $t\in[0,T]$,
 \begin{equation} 
                                      \label{weak11}  
(u_{t},\zeta):=
 \int_{\bR^d} u_t(x) \zeta(x) dx
=(u_{0},\zeta)+\int_{0}^{t}(u_{s}, L_{s} \zeta)
\,ds
+\int_{0}^{t}  (h_{s}f_{s}, \zeta) \,ds.  
\end{equation}

\begin{remark}
                                                \label{remark 5.8.3}
In light of Assumptions \ref{assumption 2.26.3} (i),
\ref{assumption 2.27.1} (ii), and \eqref{5.8.4},
the right-hand side of \eqref{weak11} makes sense for any
$u\in W$ and defines a continuous function of $t$.
Therefore, for any $h\in B((0,T),\cG)$ and $\zeta\in C^{\infty}_{0}(\bR^{d})$,
the function $(u_{t}[h],\zeta)$ is continuous on $[0,T]$.
\end{remark}

\begin{theorem}
                                   \label{theorem 2.27.1} 
Suppose  that 
 $W$, $\mathcal G$, $K$, $L$,   $u_0$, and $f$, 
described above,
satisfy Assumptions \ref{assumption 2.26.3} and 
\ref{assumption 2.27.1}.
Then,
for any  $g^{(1)},...,g^{(n)}\in B((0,T),\cG)$
and $\lambda_{1},...,
\lambda_{n} \, \geq 0$, the couple, consisting of the family of operators
 $\hat L_{t}$, such that
\begin{equation}
                                                       \label{5.16.3}
\hat L_{t}^{*}=L^{*}_{t}+\sum_{i=1}^{n}\lambda_{i}( g_{ t}^{(i)} 
-1),
\end{equation}  
where 1 stands for the operation of multiplying by one,
and $f$, is  $W$-regular. 
\end{theorem}
 
This theorem is proved in Section \ref{section 6.29.5}.

To state our second general result we need   one  more
assumption  on $W$.

\begin{assumption}
                                 \label{assumption 5.17.1}
  For any sequence $u^{k}\in W$ and a 
  bounded    function
$u=u_{t}(x)$, $(t,x)\in[0,T]\times\bR^{d}$, such that
$$
\int_{\bR^{d}}u^{k}_{t}(x)\zeta(x)\,dx\to 
\int_{\bR^{d}}u _{t}(x)\zeta(x)\,dx
$$
 for any  
  $t\in[0,T]$ and $\zeta\in C^{\infty}_{0}(\bR^{d})$,
 there exists $w\in W$ such that
$w_{t}=u_{t}$ \(a.e.\) on $\bR^{d}$ for any $t\in[0,T]$.
 \end{assumption}
 
The main consequence of Assumption \ref{assumption 5.17.1}
is the following technical result.
\begin{theorem}
                                    \label{theorem 9.17.1}
Suppose that Assumptions \ref{assumption 2.26.3} (i)
and \ref{assumption 5.17.1} are satisfied.
Let $\{L^{k}_{t},t\in(0,T)\}$, $k=0,1,...$,
be a sequence of families of linear
operators mapping $C^{\infty}_{0}(\bR^{d})$
into $B(\bR^{d})$ subject to the following conditions:

a) For each $k$, Assumption \ref{assumption 2.27.1}  \(ii\,\)
is satisfied with $L^{k}_{t}$ in place of $L _{t}$;

b) For any $\zeta\in C^{\infty}_{0}(\bR^{d})$, we have
$$
\lim_{k\to\infty}
\int_{ (0,T)\times\bR^{d} } |(L^{k}_{t}-L_{t})\zeta(x)| 
\,dt dx=0;
$$

c) For each $k=1,2,...$, there exists $u^{k}\in W$ such that
for any $\zeta\in C^{\infty}_{0}(\bR^{d})$ and $t\in[0,T]$,
 \begin{equation} 
                                      \label{9.17.1}  
(u^{k}_{t},\zeta)  
=(u_{0},\zeta)+\int_{0}^{t}(u^{k}_{s}, L^{k}_{s} \zeta)
\,ds
+\int_{0}^{t}  ( f_{s}, \zeta) \,ds.
\end{equation}

Then there exists $u^{0}\in W$ for which
\eqref{9.17.1} holds with 0 in place of $k$
for any $\zeta\in C^{\infty}_{0}(\bR^{d})$ and $t\in[0,T]$.

\end{theorem}

This theorem is proved in Section \ref{section 6.29.6}.
Theorem \ref{theorem 9.17.1} allows us to
improve the result of Theorem \ref{theorem 2.27.1}
under slightly heavier assumptions. (The conjecture is that,
actually, Assumption \ref{assumption 5.17.1} is not 
necessary in Theorem \ref{theorem 10.6.2}.)

\begin{theorem}
                                              \label{theorem 10.6.2}
Suppose  that 
 $W$, $\mathcal G$, $K$, $L$,   $u_0$, and $f$, 
described above,
satisfy Assumptions \ref{assumption 2.26.3}, 
\ref{assumption 2.27.1}, and \ref{assumption 5.17.1}.
 Let  $g^{(1)},...,g^{(n)}\in B((0,T),\cG)$ and
$\lambda_{1}(t),...,
\lambda_{n}(t)$ be nonnegative bounded measurable functions.
Then for any $h\in B((0,T),\cG)$
there exists $u\in W$ such that
\eqref{weak11} holds for any $\zeta\in C^{\infty}_{0}(\bR^{d})$ 
and $t\in[0,T]$
with 
$$
L _{s}+\sum_{i=1}^{n}\lambda_{i}(s)( g_{ s}^{(i)} 
-1),
$$
 in place of $L_{s}$.
\end{theorem}

Proof. For real variable $r$ and integer $k\geq 1$ set $\kappa_{k}(r)
=[kr]/k$, where $[r]$ stands for the integer part of $r$. 
Note that $|\kappa_{k}(r)-r|\leq1/k$, 
for any $r \in \bR$.
Set
$$
J^{0}_{t}=\sum_{i=1}^{n}\lambda_{i}(t)( g_{ t}^{(i)} 
-1),\quad J^{k}_{t}
=\sum_{i=1}^{n}\kappa_{k}(\lambda_{i}(t))( g_{ t}^{(i)}  
-1).
$$
Observe that, for an   integer  $N$, which is larger than
all $\lambda_{i}(t)$, we have 
$$
J^{k}_{t}
=\sum_{i=1}^{n}\sum_{j=1}^{Nk}(  j   /k)I_{  \{  \kappa_{k}(\lambda_{i}(t))
=j/k  \} } \, ( g_{ t}^{(i)}  
-1)=\sum_{i=1}^{n}\sum_{j=1}^{Nk}(j/k)
( g_{ t}^{(ijk)}  
-1),
$$
where $g_{ t}^{(ijk)}=g_{ t}^{(i)}$ if $\kappa_{k}(\lambda_{i}(t))
=j/k$ and $g_{ t}^{(ijk)}=1$ otherwise.

It follows by Theorem \ref{theorem 2.27.1} that for any $k\geq 1$
there exists $u^{k}\in W$ such that
\eqref{weak11} holds for any $\zeta\in C^{\infty}_{0}(\bR^{d})$ 
and $t\in[0,T]$
with 
 $
L _{s}+J^{k}_{s}
 $
 in place of $L_{s}$.

Furthermore, for any $\phi\in C^{\infty}_{0}(\bR^{d})$
$$
 \int_{ (0,T)\times\bR^{d} } |(J^{0}_{ t }- J^{k}_{ t   })\phi(x)| 
\,dt dx
$$
$$
\leq\sum_{i=1}^{n}\int_{ (0,T)\times\bR^{d} }|\lambda_{i}(t)-\kappa_{k}(\lambda_{i}(t))|
\,|(g^{(i)}_{t}-1)\phi(x)|\,dtdx
$$
$$
\leq (2T/k)\int_{  \bR^{d} }|\phi(x)|\,dx
$$
 (recall that the Jacobian of $g^{(i)}_{t}$ is one) 
 which tends to zero as $k\to\infty$. An application
of Theorem \ref{theorem 9.17.1} finishes the proof of the
present theorem.       \qed

Next, let $\frN$ be a   subset of the space
of affine 
transformations of $\bR^{d}$ and
 suppose that $\cG$ in the beginning of the section is given as
\begin{equation}
                                                     \label{5.5.1}
\cG= \{e^{t\nu}:t\in\bR,\nu\in\frN\},
\end{equation}
where by $e^{t\nu}$ we mean a transformation $
g(t)$
 defined as a unique solution
of the equation
\begin{equation}
                                           \label{2.27.5}
g(t)=1+\int_{0}^{t}\nu g(s)\,ds.
\end{equation}
  Also for any $\nu \in \frN$ we
introduce a mapping $\nu^{0}$ by the formula
$$
\nu^{0}x=\nu x-\nu0.
$$
 Notice that the
$\nu^{0}$'s are linear mappings,  which we identify with matrices
in a usual way. {\em
Of course, we keep the assumption that $\cG$ is a commutative 
group  of volume-preserving
transformations. }

 Note in passing that, 
  in case $\cG$
is given by \eqref{5.5.1}, the volume-preserving assumption
is satisfied if and only if $\tr \nu^{0}=0$ for any
$\nu\in\frN$.  
 Interestingly enough, this
``if and only if'' 
statement will never be used in the
future.  
  
 With any $\nu \in \frN$  we 
 associate an operator $M_{\nu} $ acting on  
 functions $\phi : \bR^d \to \bR$ by the formula 
$$
M_{\nu} \phi (x) =  \lim_{\varepsilon\downarrow0}
\frac{1}{\varepsilon^{2}}[\phi(e^{\varepsilon\nu}x)-2\phi(x)+\phi(e^{-\varepsilon\nu}x)],\;\; x \in \bR^d,
$$
whenever the limit on the right exists for all $x$.

Observe that if $\phi$ is twice continuously differentiable,
then
$$
M_{\nu} \phi (x) =
\frac{d^{2}}{(d\varepsilon)^{2}}\phi(e^{\varepsilon\nu}x)\big|_{\varepsilon=0}
=\frac{d }{ d\varepsilon }\big\{[D_{i}\phi](e^{\varepsilon\nu}x)
(\nu e^{\varepsilon\nu}x)^{i}\big\}\big|_{\varepsilon=0}
$$
$$
=\frac{d }{ d\varepsilon}\big\{[D_{i}\phi](e^{\varepsilon\nu}x)
(\nu^{0} e^{\varepsilon\nu}x)^{i}\big\}\big|_{\varepsilon=0}
+\frac{d }{ d\varepsilon}\big\{ [D_{i}\phi](e^{\varepsilon\nu}x)\big\}
(\nu 0)^{i} \big|_{\varepsilon=0}
$$
$$
=(\nu^{0} x)^{i}(\nu x)^{j} D_{ij}\phi(x)  
+(\nu^{0}\nu x)^{i}D_{i}\phi(x)
+(\nu 0)^{i} (\nu x)^{j} D_{ij}\phi(x)
$$
$$
 =(\nu  x)^{i}(\nu x)^{j} D_{ij}\phi(x)  
+(\nu^{2}  x-\nu0)^{i}D_{i}\phi(x).
$$

\begin{example}
                                         \label{example 2.27.1}
Let $l$ be a unit vector in $\bR^{d}$ and 
define a transformation $\nu=\nu_{l}$ by $\nu_{l} x\equiv l$
on $\bR^{d}$. Then \eqref{2.27.5} becomes
$$
g(t)x=x+\int_{0}^{t}\nu g(s)x\,ds=
x+\int_{0}^{t}l\,ds=x+tl.
$$
Observe that in this example,   
  for smooth  $\phi$, we have  
 $  M_{\nu} \phi (x)  
= D^{2}_{l}\phi(x)$.
Thus, if $\frN=\{\nu_{l}:l\in\bR^{d},|l|=1\}$, 
then $\cG$ is the set of shifts of $\bR^{d}$ and $\cG$
is a commutative group. Just in case, observe that,
for such $\frN$,   $\nu_{l_{1}}\nu_{l_{2}}\ne
\nu_{l_{2}}\nu_{l_{1}}$ unless $l_{1}=l_{2}$
although  $e^{t\nu_{1}}e^{t\nu_{2}}
=e^{t\nu_{2}}e^{t\nu_{1}}$ always.
\end{example}

\begin{example}
                                   \label{example 2.27.2}
Let $\nu x=Qx$, where $Q$ is a skew-symmetric $d\times d$-matrix.
Then  $g_t x = e^{t \nu} x = (\exp [tQ])x$, where $\exp [tQ]$ is an orthogonal matrix.
In this example, for smooth $\phi$,
$$
M_{\nu} \phi (x) 
=(Qx)^{i}(Qx)^{j}D_{ij}\phi(x)+(Q^{2}x)^{i}D_{i}\phi(x).
$$
\end{example}

\begin{theorem}
                                         \label{theorem 2.29.1}
 
 Suppose that 
$W$, $\mathcal G$, $K$, $L$, $  u_0 $ and $f$
satisfy Assumptions \ref{assumption 2.26.3} and 
\ref{assumption 2.27.1}
 with $\cG$ from \eqref{5.5.1}
and suppose that $W$ also satisfies
  Assumption  \ref{assumption 5.17.1}. 
Then,
for any  $\mu^{(1)},...,\mu^{(n)}\in B((0,T),\frN)$
equation \eqref{2.28.1} with
$$
L^{*}_{t}+\sum_{i=1}^{n}M_{\mu_{ t}^{(i)}}
$$
in place of $L^{*}_{t}$ has a solution in   $W$.  
\end{theorem}

We prove this theorem in Section \ref{section 6.29.6}.

\begin{remark} \label{sys}
  We concentrate on the case of scalar equations 
\eqref{2.28.1} only to slightly simplify the presentation.
The results similar to Theorems
\ref{theorem 2.27.1}, \ref{theorem 10.6.2}, and \ref{theorem 2.29.1}
also hold for  systems,  when $u_{t}(x)$ are vector- rather
than real-valued functions. The reader will easily
adjust our proofs to the case of systems.
 
\end{remark}

\begin{example}
                                        \label{example 2.29.4}
Let $d=2$, $\alpha\in(0,1)$, and $L_{t}=\Delta$. We know that
for any
\begin{equation}
                                            \label{3.1.1}
f\in B_{c}((0,T),C^{\infty}_{0} (\bR^{2}))  
\end{equation}
   the equation
\begin{equation}
                                                    \label{5.5.3}
u_{t}(x)=\int_{0}^{t}[\Delta u_{s}(x)+f_{s}(x)]\,ds,\quad t\leq T,
x\in\bR^{2},
\end{equation}
 has a  unique continuous solution 
such that  
$$
\sup_{(t,x) \in [0,T]\times\bR^{2}}|u_{t}(x)|+
\sup_{t\in [0,T] }\int_{\bR^{2}}|u_{t}(x)|\,dx 
$$
\begin{equation}
                                            \label{3.1.2}
\leq N_{0} \bigg[
\int_{0}^{T}\int_{\bR^{2}}|f_{t}(x)| \,dxdt+
\sup_{(t,x) \in [0,T]\times\bR^{2}}|f_{t}(x)|\bigg],
\end{equation}  
\begin{equation}
                                             \label{2.29.5}
\sup_{t \in [0,T]}[   D^{2}_{l} 
u_{t} ]_{C^{\alpha}(\bR^{2})}
\leq N_{\alpha}\sup_{t \in [0,T]}
[ f_{t} ]_{C^{\alpha}(\bR^{2})} 
\end{equation}
 for any   $l\in S_{1}=\{|x|=1\}$,  where   
     $N_{0}$ and $N_{\alpha}$ are some constants.

We   claim  that, if  \eqref{3.1.1} holds, 
the equation
$$
u_{t}(x)=\int_{0}^{t}[\Delta u_{s}(x)+M u_{s}(x)+f_{s}(x)]\,ds,
$$
where
$$
M\phi(x)=(x^{2})^{2}D_{11}\phi(x)-2x^{1}x^{2}D_{12}
\phi(x)+(x^{1})^{2}D_{22}\phi(x)
$$
$$
 - x^{1}D_{1}\phi(x) - x^{2}D_{2}\phi(x),
$$
has a continuous solution, which satisfies estimates \eqref{3.1.2}
and \eqref{2.29.5}
(with the same $N_{0}$ and $N_{\alpha}$).

  With the goal of applying Theorem
\ref{theorem 2.29.1},  fix $f$ as in \eqref{3.1.1} and
 denote by $A_{0}$ and $A_{\alpha}$
the right-hand sides of \eqref{3.1.2} and \eqref{2.29.5},
respectively. Then introduce
 $$
W=\{u\in { B([0,T]\times\bR^{2}):} 
 \; u_t \in C^{2+ \alpha}(\bR^d), \; t \in [0,T],   
  \sup_{(t,x)\in [0,T]\times\bR^{2}}|u_{t}(x)|  
$$
$$+
\sup_{t\in [0,T] }\int_{\bR^{2}}|u_{t}(x)|\,dx \leq A_{0},
 \sup_{t \in [0,T]}[ D^{2}_{l} u_{t}
]_{C^{\alpha}(\bR^{2})}
\leq A_{\alpha}\,\,
 \forall l\in S_{1} \},
$$
and let $\frN=\{t Q :t\in\bR\}$, where 
$Q = (Q_{ij})$ is a $2 \times 2$-matrix, $Q^{ii}=0$,
$Q^{12}=1$,
$Q^{21}=-1$,  $i=1,2$.  
Note that since $Q$ is skew-symmetric, $\cG=\{e^{tQ};t\in \bR\}$
is the group of rotations of $\bR^{2}$  about
the origin. 

 In light of Example \ref{example 2.27.2} 
and Theorem
\ref{theorem 2.29.1}, to prove our
claim, it suffices to check that 
Assumptions \ref{assumption 2.26.3}, 
\ref{assumption 2.27.1}, and
\ref{assumption 5.17.1} 
  are satisfied
for the above $W$ and $\frN$, $u_{0}=0$
and $\Delta$ in place of $L_{t}$.
 
 Assumption  \ref{assumption 2.26.3} (i) is obviously
satisfied. Assumption  \ref{assumption 2.26.3} (ii)   is  
satisfied since, for instance,
$$
\sup_{t \in [0,T]}[ D^{2}_{l} 
Eu_{t} ]_{C^{\alpha}(\bR^{2})}
\leq \sup_{t \in [0,T]}E[ D^{2}_{l} 
u_{t}]_{C^{\alpha}(\bR^{2})}\leq A_{\alpha}.
$$ 
 Moreover, using that   $|g x| = |x|$, $g\in \cG$,  we deduce 
that for  any bounded  measurable $\cG$-valued function
$g_{t}$  given on $[0,T]$
$$
\big[ D^{2}_{l}( u_{t}(  g_t \cdot)) 
 \big]_{C^{\alpha}(\bR^{2})}
=\big[ (D^{2}_{g_{t}l} u_{t})(  g_t \cdot) 
 \big]_{C^{\alpha}(\bR^{2})}
=\big[  D^{2}_{g_{t}l} u_{t}  
 \big]_{C^{\alpha}(\bR^{2})}\leq A_{\alpha}.
$$
 
By adding to this that
$$
 \int_{\bR^{2}}|u_{t}(g_tx)|\,dx 
= \int_{\bR^{2}}|u_{t}( x)|\,dx 
$$
since $\det g_{t}=1$, we conclude that the function
 $u_{t} ( g_{t}x)$ is  in $ W$ and 
Assumption \ref{assumption 2.26.3}
is satisfied. 

  Assumption \ref{assumption 2.27.1} (ii) is obviously
 satisfied and requirement
(i) is satisfied since the Laplacian is rotation invariant.
As long as Assumption \ref{assumption 2.27.1} (iii)
is concerned, observe that, for any  $h  \in B((0,T) , \cG )$,
we have
 ${h }_t f_t =f_{t}({h }_{t}\cdot) \in  
B_{c}((0,T),C^{\infty}_{0} (\bR^{2}))$, so that
equation \eqref{5.5.3} with ${h }_s f_s$ in place of $f_{s}$
 has a  unique continuous solution and estimates
\eqref{3.1.2} and \eqref{2.29.5} are valid with
${h }_t f_t$ in place of $f_{t}$. As is seen from the above
arguments, this replacement
does not change the right-hand sides of \eqref{3.1.2} and \eqref{2.29.5},
which implies that Assumption \ref{assumption 2.27.1} (iii)
is satisfied. That Assumption \ref{assumption 2.27.1} (iv)
is satisfied is a simple consequence of the maximum principle.

  To check Assumption \ref{assumption 5.17.1},    
 we consider a sequence $u^k$ which converges in the 
specified  weak
sense to a function $u$ defined on $[0,T] \times \bR^2$.
We fix $t \in [0,T]$.
Possibly passing to a subsequence and using  the 
Arzel\`a-Ascoli
theorem, we find  that there exists $w_t  \in C^{2+
\alpha}(\bR^2)$ such that, along the   subsequence, 
$u^k_t $, $D_i u^k_t $, and $D_{ij} u^k_t$  converge to $w_t $,
$D_i w_t$, and  $D_{ij} w_t $, respectively, uniformly on each compact 
subset
of $\bR^2$.  In principle it could happen 
that along a different subsequence 
$u^k_t $, $D_i u^k_t $, and $D_{ij} u^k_t$  converge to $w'_t $,
$D_i w'_t$, and  $D_{ij} w'_t $ uniformly on each compact subset
of $\bR^2$ and $w_{t}\ne w'_{t}$. However,
along both subsequences
$$
\int_{\bR^{d}}u^{k}_{t}\zeta\,dx\to\int_{\bR^{d}}u _{t}\zeta\,dx
$$
for any   
 $\zeta\in C^{\infty}_{0}(\bR^{d})$.
It follows that 
$$
w_{t}= w'_{t}=u_{t}
$$
 in $\bR^{d}$ almost everywhere,
and, since $w_{t}$ and $w'_{t}$ are continuous, $w_{t}= w'_{t}$
 everywhere.

Thus, for each $t \in [0,T]$,
the sequences
$u^k_t $, $D_i u^k_t $, and $D_{ij} u^k_t$  converge to $w _t $,
$D_i w _t$, and  $D_{ij} w _t $, respectively, uniformly on 
each compact subset
of $\bR^2$ as $k\to\infty$. Since $u^k_t (x)$ are 
Borel measurable  as functions of $(t,x)$,
so is $w_{t}(x)$. The fact that $w$ satisfies the inequalities
entering the definition of $W$  is obvious.
 This proves our claim.
\end{example}

\begin{remark}
                                        \label{first}
In Theorem \ref{theorem 2.29.1} we could consider more 
general operators like 
\begin{equation}
                                              \label{d22}
L^{*}_{t}+\sum_{i=1}^{n}M_{\mu_{ t}^{(i)}} +
 \sum_{j=1}^{m} F_{\nu_{ t}^{(j)}},
\end{equation} 
where  $F_{\nu_{ t}^{(j)}}$ are 
 first-order operators defined by
$$
F_{\nu} \phi (x) 
= [D_{i}\phi](x)
(\nu x)^{i}.
$$
 The conclusion of Theorem \ref{theorem 2.29.1} 
remains true
since the   substitution  $v_{t}(x)=u_{t}(g^{(1)}(t)
\cdot...\cdot g^{(m)}(t)x)$,
where $g^{(i)}(t)=e^{\nu^{(i)}t}$, converts the equation for $v_{t}(x)$
not   containing  $F$'s into an equation for $u_{t}(x)$
with the additional first-order terms. Of course, the free
term will change. But it will satisfy the same estimates
as before the above change of variables. 
 \end{remark}
 
\begin{example}
                                              \label{example 5.6.1}
  As mentioned in Remark  \ref{sys}  results similar to Theorems
\ref{theorem 2.27.1}, \ref{theorem 10.6.2}
 and \ref{theorem 2.29.1}
also hold for  systems. 
 Without going into too much detail, 
we just give an example
of the following hyperbolic system in $\bR^{2}$:
\begin{equation}
                                                     \label{5.7.1}
\partial_{t}w_{t}(x)= v_{t}(x),\quad
\partial_{t}v_{t}(x)=D _{11} w_{t}(x)
\end{equation}
on $[0,T]\times \bR^{2}$ with initial condition 
$w_{0}(x)=\zeta(x^{1})\eta(x^{2}),
v_{0}(x)= \zeta'(x^{1})\eta(x^{2})$,
where $\zeta,\eta\in C^{\infty}_{0}(\bR)$ are fixed function (of one variable
and $\zeta'$ is the derivative of $\zeta$).
Assume that $\zeta,\eta\geq0$. Of course, $x^{2}$ enters system \eqref{5.7.1}
only as parameter.

Take $\frN$ and $\cG$ from Example \ref{example 2.27.1} and
define $W$
 as the collection of 
Borel $\bR^{2}$-valued functions $u_{t}(x)=(\psi_{t}(x),\phi_{t}(x))$  
on $[0,T]\times\bR^{2}$ 
such that 
$$
\psi\geq 0\quad \text{in}\quad [0,T]\times\bR^{2}\,\,\text{(a.e.)},
$$ 
$$
\int_{[0,T]\times\bR^{2} } \psi_{t}(x) \,dxdt
\leq T\int_{    \bR^2   }w_{0}(x) \,dx,
$$
$$
\int_{[0,T]\times   \bR^2}   |\phi_{t}(x)| \,dxdt
\leq T\int_{  \bR^2}  |v_{0}(x)| \,dx.
$$
Of course, given an $\bR^{2}$-valued function  $(\psi(x),\phi (x))$   
and $g\in\cG$, we define $g(\psi (x),\phi(x))=(\psi (gx),\phi (gx))$.
Then, obviously, Assumption \ref{assumption 2.26.3}
is satisfied. Also observe that since by definition
$g(\psi(x),\phi(x))=(\psi(gx),\phi (gx))$, the operator 
$M_{\nu_{l}}$ from Example \ref{example 2.27.1} will act 
on vector-valued functions by the formula  $M_{\nu_{l}}
 (\psi(x),\phi(x))=(D^{2}_{l}\psi(x),D^{2}_{l} \phi(x))$
if $\psi$ and $\phi$ are smooth enough.

Next, we define $L_{t}$ to be a $2\times 2$ matrix whose entries
are operators: $L_{t}^{11}=L_{t}^{22}=0$, $L_{t}^{21}=D_{11}$,
and $L_{t}^{12}$ is a unit operator. Finally, set $f\equiv0$.

Then system \eqref{5.7.1} in the integral form
becomes \eqref{2.28.1} and, for any 
bounded measurable
$\cG$-valued functions ${h}={h}_{t}$, $t\in(0,T)$, it 
has a solution 
\begin{equation}
                                              \label{6.30.5}
u_{t}(x)=(\zeta(x^{1}+t)\eta(x^{2}), 
\zeta'(x^{1}+t)\eta(x^{2}))
\end{equation}
(independent of ${h}$). This shows that Assumption \ref{assumption 2.27.1}
is also satisfied. Assumption  \ref{assumption 5.17.1} 
  is   easily
verified as well, 
 and by a vector-valued counterpart of Theorem \ref{theorem
2.29.1} we obtain that the  parabolic system
$$
\partial_{t}w_{t}(x)=v_{t}(x)+\Delta w_{t}(x),
\quad \partial_{t}v_{t}(x)= D _{11} w_{t}(x)+\Delta v_{t}(x),
$$
$t\in[0,T]$, $x\in\bR^{2}$, with initial data
$w_{0}(x)=\zeta(x^{1})\eta(x^{2})$ and
$v_{0}(x)= \zeta'(x^{1})\eta(x^{2})$ has a solution (in the sense explained
after Assumption \ref{assumption 2.27.1}) belonging to $W$.

 In particular, for this solution
$w_{t}(x)\geq0$ (a.e.). Actually, this result comes as no surprise
since $(w_{t},v_{t})=T_{t}u_{t}$, where $u_{t}$
is defined in \eqref{6.30.5} and $T_{t}$ is the heat semigroup
 acting on $\bR^2$-valued functions.
 We just wanted to show
that our main results are applicable to systems of equations
and not only in what concerns a priori 
estimates  for scalar equations.
\end{example}

\begin{example}
                                      \label{example 10.6.1}
Consider the following hyperbolic system taken from  
\S 7.3.3 of \cite{Ev}
\begin{equation}
                                             \label{10.6.1}
\partial_{t}u^{ r}_{t}(x)+B^{rk}_{j}D_{j}u^{k}_{t}(x) 
=g^{r}_{t}(x)
\end{equation}
$r=1,...,m$, in $(0,T)\times\bR^{d}$ with 
zero initial condition, where the $m\times m$ constant
matrices $B_{j}:=(B^{rk}_{j})$, $j=1,...,d$,
 are such that for any
$\xi\in\bR^{d}$, the matrix $\xi^{j}B_{j}$ has $m$ 
{\em real\/} eigenvalues. Assume that $g_{t}(x)
=(g^{r}_{t}(x))$ 
is an $\bR^{m}$-valued measurable functions such that 
$$
\int_{0}^{T}\|g_{t}\|_{H^{s} (\bR^{d}; \bR^m)}^{2}\,dt=A<\infty,
$$
where $s>m+d/2$ and $H^{s} (\bR^{d}; \bR^m)=W^{s}_{2}(\bR^{d}; \bR^m)$ are the usual fractional
Sobolev spaces  of $\bR^m$-valued functions (see their definitions, for instance, in 
\S 5.8.4 of \cite{Ev}). By closely following the proof
of Theorem 5 in \S 7.3.3 of \cite{Ev} (given there
for $g=0$ but with nonzero initial value) one arrives 
at the conclusion that \eqref{10.6.1} with zero initial condition
has a unique solution in class $W$, which consists
of measurable functions $u=u_{t}(x)$ on $[0,T]\times\bR^{d}$,
such that $u_{t}\in C^{0,1}(\bR^{d}; \bR^m)$  
(here $C^{0,1}(\bR^{d}; \bR^m)$ is the usual space of $\bR^m$-valued Lipschitz functions on $\bR^d$) for any $t\in[0,T]$ and
\begin{equation}
                                              \label{10.7.1}
\|u\|_{L_{2}([0,T]\times\bR^{d}; \bR^m)}
+\sup_{t\in[0,T]}\|u_{t}\|_{C^{0,1}(\bR^{d}; \bR^m)}\leq N'A,
\end{equation}
where $N'$ is a constant independent of $g$. As in Example
\ref{example 2.29.4} one checks that
the assumptions of Theorem \ref{theorem 2.29.1}
are satisfied with obvious matrix-valued first-order
differential operator and $\cG$ being the group of translations.
               
Now take a bounded measurable $d\times d$-matrix valued function
$a=a_{t}$ which is symmetric and nonnegative for any $t\in[0,T]$.
Define $\sigma_{t}=a^{1/2}_{t}$. One knows that $\sigma_{t}$
is also measurable and if $\sigma^{(i)}_{t}$ is the $i$th column of
$\sigma(t)$, $i=1,..,d$, then for smooth $\phi=\phi(x)$
$$
a^{ij}_{t}D_{ij}\phi=\sum_{i=1}^{d}D^{2}_{\sigma_{t}^{(i)}}\phi
$$
(cf. Example \ref{example 2.27.1}). Therefore,
by Theorem \ref{theorem 2.29.1} system \eqref{10.6.1}
with the additional terms on the right-hand side $a^{ij}_{t}D_{ij}u^{r}_{t}(x)$
has a solution of class $W$. In particular, estimate \eqref{10.7.1}
holds for the solution of the new system with the {\em
same\/} right-hand side. Observe that
the system is of unknown type,
because  no nondegeneracy assumption is imposed on $a_{t}$.

It is worth mentioning that  the fact that
estimate \eqref{10.7.1} holds for the new system with a constant $N'$
independent of $a$
 can also be obtained by
 closely following the proof
of Theorem 5 in \S 7.3.3 of \cite{Ev}.

\end{example}

\begin{remark}
\label{remark 10.7.1}

It could be that in each of the above examples one
can prove our assertions by examining the classical proofs.
However, the whole point is that under some
easily verified conditions we have
a unified method of adding new term into
the equations without caring much as of why an how the sets $W$
were proved to be appropriate in any particular problem.

 Just in case, we recall that all equations
are understood in a weak sense as in \eqref{weak11}.
\end{remark}

\mysection{Proof of Theorem \protect\ref{theorem 2.27.1}}
                                              \label{section 6.29.5}

We need some preparations. Again take  
  independent and identically exponentially distributed  
 with parameter  $\lambda> 0$  random variables
$\tau_{1},\tau_{2},...$ defined on a 
 probability space $(\Omega,\cF,P)$ and construct $\pi_{t}$
as in Section \ref{section 1.25.1}. 
For $t\geq0$ introduce $\cF_{t}$ as the smallest $\sigma$-fields
in $\Omega$  containing all sets
of the form $\{\omega:\pi_{s}(\omega)=k\}$, $s\leq t$, $k=0,1,...$.
Since, for $t>s$, $\pi_{t}-\pi_{s}$ is independent
of $\pi_{r}$, $r\leq s$, $\pi_{t}-\pi_{s}$ and $\cF_{s}$
are independent.

Also take
$g \in B((0,T),\cG)$, extend it to $[0,\infty)$
by setting $g_{0}=1$ and $g_{t}=1$ for $t\geq T$,
where $1$ is the operator of multiplying by 1,
 and define  $h_{t}=h_{t}(\omega) \in \cG$ 
for $t\geq 0$ and $\omega\in \Omega$
by
\begin{equation}                                             \label{h1}
h_{t}= g_{\sigma_{n}} h_{\sigma_{n}-}   
\quad\text{for}\quad t\in
[\sigma_{n} ,\sigma_{n+1} ),
\end{equation}
$n=0,1,...$, where $\sigma_{0 }-=0-:=0$ and $h_{0}x 
 :\equiv x$, $x\in\bR^{d}$.
In other terms,
$$
h_{t}= \prod_{n\leq \pi_{t}}g_{\sigma_{n}}=
\prod_{n\leq \pi_{t}}g_{\sigma_{n}\wedge t}.
$$
 Observe that
the random variables
$\sigma_{n}\wedge t$ are $\cF_{t}$-measurable because,
for constant $c\geq0$,
the set $\{\omega:\sigma_{n}(\omega)\wedge t\leq c\}$
coincides with $\Omega$ if $c\geq t$, and if $c\in[0,t)$,
this set is $\{\omega:\sigma_{n}(\omega) \leq c\}
=\{\omega:\pi_{c}(\omega) \geq n\}\in \cF_{c}\subset\cF_{t}$.
Since $g_{t}$ is measurable, $g_{\sigma_{n}\wedge t}$
is $\cF_{t}$-measurable.
It follows that $h_{t}$ is $\cF_{t}$-measurable for each $t$,
or, in other words, the process $h_{t}$ is 
$\cF_{t}$-adapted.

 The construction of the stochastic process
$h$ with values in $\mathcal G$ is inspired by the one of the simpler  process $b_t$
used in the proof of Lemma 
\ref{lemma 1.25.4}.

   Also note that the number of jumps of $\pi_{t}$
on $[0,T]$ is finite and, therefore, $h_{t} (\omega) $ is bounded on $[0,T]$
for any $\omega$.

 Before the next result recall  that the notation  $u_t[h]$
is introduced in Assumption \ref{assumption 2.27.1},
 and $(u_{t},\zeta)$  
in \eqref{weak11} and, according to what
is said in the beginning of Section \ref{section 6.29.3},
  $gx$ is the image of $x$ under mapping $g\in\cG$.

\begin{lemma}
                                               \label{lemma 2.28.3}
 Let $h$   be   introduced
by  \eqref{h1} and let $\hat h\in B((0,T),\cG)$. Then

 \(i\,\) For any $\zeta\in C^{\infty}_{0}(\bR^{d})$, the process
$\eta_{t}:=(u_{t}[h(\omega)\hat h],\zeta)$, $t\in[0,T]$ 
is continuous and $\cF_{t}$-adapted.

\(ii\,\) For any nonrandom bounded measurable $\cG$-valued function $\beta_{t}$,
$t\in (0,T]$, the function
\begin{equation}
                                                  \label{5.15.1}
u_{t}[h(\omega)\hat h](h^{-1}_{t-}(\omega)\beta_{t}x)
\end{equation}
 is 
$\cF_{T}\times \cB([0,T]\times\bR^{d})$-measurable and
belongs to $W$ for any $\omega$.

\end{lemma}

Proof.   We will see that the assertions of the lemma
hold true no matter which $f$, satisfying \eqref{5.8.4},
is taken in \eqref{2.28.1} in construction of $u_t[h]$.
 Therefore, by replacing
$f$ in \eqref{2.28.1} with $f'=\hat h f$ we reduce the
 general situation to the one where $\hat h\equiv1$,
which we assume henceforth.

 (i). The continuity of $\eta_{t}$
 follows from Remark \ref{remark 5.8.3}. 
 To investigate its measurability
properties, we need
the separable
 Banach space $L_{1}((0,T),\cG)$
of measurable   and integrable  $\cG$-valued functions on $[0,T]$.
Notice that any element $\alpha\in\cG$ is an affine transformation
and $\alpha x$ has a unique representation as $a_{\alpha}x+b_{\alpha}$,
where $a_{\alpha}$ is a linear mapping and $b_{\alpha}$ is a vector.
The norms of $a_{\alpha}$ and $b_{\alpha}$ are well defined and we make
the space, say   $\Lambda$,
  of affine transformation of $\bR^{d}$ a linear
normed space by setting
$$
|\alpha'-\alpha''|=|a_{\alpha'}-a_{\alpha''}|+|b_{\alpha'}-b_{\alpha''}|.
$$
After that we introduce the norm in the linear
space $L_{1}((0,T),\cG)$
by setting
$$
\|\alpha\|_{L_{1}((0,T),\cG)}=
\int_{0}^{T}|\alpha_{t}|\,dt.
$$
As any $L_{1}$-space relative to Lebesgue
measure of functions on $(0,T)$ with values
in {\em finite-dimensional\/} spaces,    
the space $L_{1}((0,T),\cG)$ is  Polish.

Next, we take continuous   $\Lambda$-valued functions
$\phi^{{m}}(\alpha)$,
${m}=1,2,...$, on  $\Lambda$  each of which is bounded
and such that 
$\phi^{{m}}(\alpha)=\alpha$
for $|\alpha|\leq {m}$.

Observe  that, if $\alpha^{n} \in L_{1}((0,T),\cG)$, $n=0,1,...$,
are such that $\alpha^{n} \to \alpha^{0} $  in
$L_{1}((0,T),\cG)$ as $n\to\infty$, then, for any fixed
$m=1,2,...$, $\alpha^{mn}:=\phi^{m}(\alpha^{n}) \in
B((0,T),\cG)$, so that $u_{t}[\alpha^{mn}](x)$
are well defined. We claim that in this situation
$u_{t}[\alpha^{mn}](x)
\to u_{t}[\alpha^{m0}](x)$ uniformly on $[0,T]\times\bR^{d}$
as $n\to\infty$.

To prove this claim,
thanks to  Assumption \ref{assumption 2.27.1}  (iv),
it suffices to show that, for any fixed $m$,
\begin{equation}
                                              \label{8.23.1}
I_{n}:=\int_{0}^{T}
\sup_{y \in \bR^{d}}|f_{r}( \alpha^{mn}_{r}y)
-f_{r}( \alpha^{m0} _{r}y)| \,dr
\to0
\end{equation}
as $n\to\infty$. As usual, it suffices to prove
\eqref{8.23.1} assuming that $\alpha^{n}_{t}\to \alpha^{0}_{t}$
for almost any $t$. For such $t$ 
and any $y$ we have $f_{t}( \alpha^{mn}_{t}y)-
f_{t}( \alpha^{m0} _{t}y)\to0$
by continuity.   Furthermore  the functions
$f_{t}( \alpha^{mn}_{t}y)-f_{t}( \alpha^{m0} _{t}y)$ are supported
in the same ball
and are uniformly continuous ($t$ and $m$ are fixed). Therefore, the
convergence $f_{t}( \alpha^{mn}_{t}y)-f_{t}( \alpha^{m0} _{t}y)\to0$
is uniform on $\bR^{d}$, and this implies \eqref{8.23.1}
by the dominated convergence theorem.

Hence, $u_{t}[\phi^{m}(\alpha)](x)$ is   continuous  with respect to
  $\alpha\in L_{1}((0,T),\cG)$ uniformly with respect to $(t,x)$.

Next, coming back to $h(\omega)$ observe that
for any $\alpha\in L_{1}((0,T),\cG)$ the   random  function
$$
\rho(\alpha,h   ):=\int_{0}^{T}|\alpha_{t}-h_{t} |\,dt =\sum_{n\leq\pi_{T}}\int_{\sigma_{n}{\wedge T}}^{
 \sigma_{n+1} \wedge T}
\Big|\alpha_{t}-\prod_{i\leq n}g_{\sigma_{i}{\wedge T}}\Big|\,dt
$$
is  $\cF_{T}$-measurable. Therefore, we have
$$
\{\omega:\rho(\alpha,h(\omega))\leq\rho\}\in\cF_{T}
$$
for any $\alpha\in L_{1}((0,T),\cG)$ and $\rho>0$,
Since $L_{1}((0,T),\cG)$ is a Polish space, we get that
$h(\omega)$ is an $\cF_{T}$-measurable $L_{1}((0,T),\cG)$-valued
function.

Now   we conclude that,
since $u_{t}[\phi^{m}(\alpha)](x)$ is   continuous in $\alpha$
and $h(\omega)$ is $\cF_{T}$-measurable,
$u_{t}[\phi^{m}(h(\omega))](x)$ is $\cF_{T}$-measurable.
By observing that $h_{t}(\omega)$ is bounded for each
$\omega$ by definition, we conclude that
$u_{t}[\phi^{m}(h(\omega))](x)\to u_{t}[h(\omega)](x)$
as   $m \to \infty$ uniformly with respect to $(t,x)$
in the sense that, actually, for each $\omega$, there is
$n(\omega)$ such that 
$u_{t}[\phi^{m}(h(\omega))](x)= u_{t}[h(\omega)](x)$
on $[0,T]\times\bR^{d}$,   for any $m \ge n(\omega)$. Anyhow, this proves
that $u_{t}[h(\omega)](x)$ is $\cF_{T}$-measurable.

By fixing $t\in[0,T]$, replacing $T$ above with $t$,
and applying the above argument to 
$(u_{t}[h(\omega) ],\zeta)$, we get that $\eta_{t}$
is  $\cF_{t}$-measurable. This proves (i).

(ii). By the above
  $u_{t}[\phi^{m}(\alpha)](x)$ is   continuous in $\alpha
 \in L_{1}((0,T),\cG)$ for any $m$ and,
 by definition, it is Borel measurable
with respect to $(t,x)$, for any fixed $\alpha$.
A general simple result then tells us that
$u_{t}[\phi^{m}(\alpha)](x)$ is Borel measurable
in $(\alpha,t,x)$, that is ${\mathcal
B}(L_{1}((0,T),\cG)) \times  [0,T] \times
\bR^d)$-measurable.  The mapping $(\omega,t,x)
\to (h(\omega),t,x)$ is also measurable,
and since the superposition of a Borel measurable 
and a measurable function is measurable,
$u_{t}[\phi^{m}(h(\omega))](x)$ is $\cF_{T} \times 
{\mathcal B}([0,T] \times \bR^d)  $-measurable with respect
to $(\omega,t,x) $. By letting $m\to\infty$, we conclude that
  $u_{t}[h(\omega)](x)$ $\cF_{T} \times 
{\mathcal B}([0,T] \times \bR^d)  $-measurable
 with respect to $(\omega,t,x) $.

Next, $h^{-1}_{t-} (\omega)
\beta_{t}x$ is  a measurable function of $(\omega,t,x)$,
and by the properties of superpositions of measurable functions
$$
u_{t}[  h(\omega) ]
(h^{-1}_{t-}(\omega) \beta_{t}x )
$$
is measurable with respect to $(\omega,t,x ) $,
that is, it
is $\cF_{T}\times \cB([0,T]\times\bR^{d}) $-measurable indeed.

The fact that it belongs to $W$ for each $\omega$
follows directly from Assumptions 
\ref{assumption 2.26.3} (iii) and \ref{assumption 2.27.1} (iii).

The lemma is proved.  \qed 
 
 Next, we need the notion of predictable $\sigma$-field $\cP$.
This is the smallest $\sigma$-field of subsets
of $\Omega\times(0,\infty)$ containing all sets
of the form $B=A\times(s,t]$, where $A\in\cF_{s}$ and
$0\leq s<t <\infty $ are arbitrary. $\cP$-measurable
functions on $\Omega\times(0,\infty)$ are called
predictable processes.  
It is convenient to speak about predictable processes given
only on $(0,T]$, we then just continue them as their
values at $T$ after that time.  It is a well-known
and easy fact that  all real-valued left-continuous,
$\cF_{t}$-adapted processes $\xi_{t}(\omega)$ given
on $\Omega\times(0,\infty)$ are predictable.
  In particular, the process $\eta_{t}$
from Lemma \ref{lemma 2.28.3} is predictable for any $\zeta$.

  A trivial example 
of predictable function is given by any (nonrandom 
Borel) measurable function on $(0, \infty) $.
It is predictable, because the smallest $\sigma$-field
containing all intervals $(s,t]$   is the Borel $\sigma$-field 
of   $(0, \infty) $.

\begin{remark}
                                               \label{remark 5.15.2}
If $\beta_{t}$ is a predictable $\cG$-valued process,
$\zeta\in C^{\infty}_{0}(\bR^{d})$, $h$ is taken from \eqref{h1},
and $\hat h\in B((0,T),\cG)$,
then $(u_{t}[h(\omega)\hat h],\zeta(\beta_{t}(\omega)\cdot)$
is predictable. This follows from the fact that 
$(u_{t}[h(\omega)\hat h],\zeta(\beta \cdot)$ is predictable
for any $\beta\in\cG$ and is continuous with respect to
$\beta$ so that it is jointly measurable with 
respect to $(\omega,t,\beta)$.

\end{remark}

We are going to use the following.

\begin{lemma}
                                        \label{lemma 2.29.2}
Let $\xi_{t}$ be a   predictable process such that
$$
E\int_{0}^{t}|\xi_{s}|\,ds<\infty.
$$
 Then
\begin{equation}
                                                 \label{3.1.11}
E\int_{(0,t]}\xi_{s}\,d\pi_{s}=
\lambda E\int_{0}^{t}\xi_{s}\,ds.
\end{equation}
\end{lemma}

   This lemma 
 follows from Theorem 16  and
the comments after it in Section III.5 on page 118 of \cite{Pr}. Since going
through the material before that theorem can be somewhat painful for
inexperienced reader we give a  
 short proof.

First of all we note that it suffices to concentrate
on bounded processes $\xi_{t}$. This follows from
the monotone convergence theorem by a routine argument.
In that case the lemma
is just Exercise 2.7.8 of \cite{Kr02} and its solution,
 given below,  is outlined in the hint to this exercise.

  One fixes $t>0$ and   introduces  two measures
on $\cF\times\cB(0,t]$
$$
\mu(B)=E\int_{(0,t]}I_{B}(\omega,r)\,d\pi_{r},\;\;
\nu(B)=\lambda E\int_{(0,t]}I_{B}(\omega,r)\,d r.
$$
When $B=A\times(a,b]$,   with $A \in \cF_a$, 
  and $a<b\leq t$, we have  $\mu(B)=
EI_{A}(\pi_{b}-\pi_{a})=\lambda P(A)(b-a)=\nu(B)$
 because   $I_A$   and $\pi_{b}-\pi_{a}$
 are independent.  The equality 
$\mu(B)=\nu(B)$ for $B=A\times(a,b]$ is also easily verified
for other dispositions of $a,b,t$.
 Thus, $\mu=\nu$ on such sets $B$.
Since the collection of such $B$ is  a $\Pi$-system
  (see the definition
of $\Pi$-system in \cite{Kr02}), by a very general
 fact from measure theory (see Lemma 2.3.18 in \cite{Kr02})
$\mu(B)=\nu(B)$ on the smallest $\sigma$-field containing all
such $B$, that is, on    $\cP$.

We thus have proved \eqref{3.1.11}
if $\xi_{r}(\omega)$ is the indicator of a predictable set.
The same equality is true if $\xi_{r}(\omega)$ is a finite linear
combination of the indicators of predictable sets
with nonrandom coefficients.
Since bounded measurable functions admit
uniform approximations by finite linear combinations
of the indicators of measurable sets, \eqref{3.1.11}
holds for all bounded predictable processes
and the lemma is proved. \qed

\begin{remark}
The reader may feel uncomfortable encountering
the above measure-theoretic arguments which we easily
avoided in Sections \ref{section 1.25.1} and \ref{section 6.29.2}.
Unfortunately, these arguments are necessary in the general 
 theory.
To see this, observe that
$$
\int_{(0,t]}\pi_{s}\,d\pi_{s}=
\sum_{n=0}^{\pi_{t}}n=(1/2)\pi_{t}(\pi_{t}+1),\quad
E\int_{(0,t]}\pi_{s}\,d\pi_{s}=  \lambda t+\lambda^{2}t^{2}/2.  
$$
At the same time
$$
E\int_{0}^{t}\pi_{s}\,ds=\int_{0}^{t}E\pi_{s}\,ds=
\lambda\int_{0}^{t}s\,ds=\lambda t^{2}/2,
$$
and \eqref{3.1.11} does not hold for $\xi_{s}=\pi_{s}$.

However,
$$
\int_{(0,t]}\pi_{s-}\,d\pi_{s}=\int_{(0,t]}[\pi_{s}-1]\,d\pi_{s}
=(1/2)\pi_{t}(\pi_{t}-1),\quad
E\int_{(0,t]}\pi_{s-}\,d\pi_{s}=\lambda t^{2}/2.
$$
By the way, one of consequences
of these calculations and Lemma \ref{lemma 2.29.2}
is that the process $\pi_{t}$ is not predictable,
although $\pi_{t-}$ is.
\end{remark}

{\bf Proof of Theorem  \ref{theorem 2.27.1}}.   
Obvious induction on $n$ allows us to concentrate on $n=1$
and assume that $\lambda=\lambda_{1}>0$.  
Next, the
requirements (i) and (ii) of Assumption \ref{assumption 2.27.1}
are obviously satisfied for the operators 
whose formal adjoints are defined in  \eqref{5.16.3}.
To check the remaining requirements, take $g\in B((0,T),\cG)$,
take $h$ and $\zeta$
as in Lemma \ref{lemma 2.28.3}, take any $\hat h\in B((0,T),\cG)$,
 and consider the process
$$
\xi_{t}=(u_{t}[h\hat h],\zeta(h_{t}\cdot)),
$$
where and below we drop the argument $\omega$ as usual.
  This process is well-defined since  changing variables 
(recall Remark \ref{remark 8.24.1}) we get
$$ |(u_{t}[h(\omega)\hat h],\zeta(h_{t}\cdot))| \le \sup_{(t,x) \in
[0,T] \times \bR^d} |u_{t}[h(\omega)\hat h](x)| \int_{\bR^d}
|\zeta(y)| dy.
$$
By the same reason the processes 
$\xi_{t,r}=(u_{t}[h\hat h],\zeta(h_{r})\cdot)$ are well defined
for $t,r\in[0,T]$. In addition, for any fixed $r\in[0,T]$,
 viewing $\zeta(h_{r}x)$
just as another $C^{\infty}_{0}(\bR^{d})$-function, for
$t\in[\sigma_{n},\sigma_{n+1})$ we obtain
$$
\xi_{t,r}=\xi_{\sigma_{n},r}+\int_{\sigma_{n}}^{t}
(u_{s}[h\hat h],L_{s}\zeta(h_{r}\cdot))
 \,ds+
\int_{\sigma_{n}}^{t}(h_{s}\hat h_{s}f_{s},\zeta(h_{r}\cdot))
 \, ds. 
$$
We substitute here $r=\sigma_{n}$ and observe that for
$t\in[\sigma_{n},\sigma_{n+1})$ the function
$h_{t}$ does not change and equals $h_{\sigma_{n}}$ and
$$
\xi_{t,\sigma_{n}}=(u_{t}[h\hat h],
\zeta(   h_{\sigma_{n}}\cdot)   )
=(u_{t}[h\hat h],   \zeta(h_{t} \cdot ) )=\xi_{t}.
$$
Then we conclude that similarly to \eqref{1.23.05},
for
$t\in[\sigma_{n},\sigma_{n+1})$,
$$
\xi_{t}=\xi_{\sigma_{n}}+\int_{\sigma_{n}}^{t}
(u_{s}[h\hat h],L_{s}\zeta(h_{s}\cdot))
 \,ds+
\int_{\sigma_{n}}^{t}(h_{s}\hat h_{s}f_{s},\zeta(h_{s}\cdot))
\,ds.
$$
At time $t=\sigma_{n+1}$ the process $h_{t}$ jumps
from    $h_{\sigma_{n+1}-}$ to 
$h_{\sigma_{n+1}}= g_{\sigma_{n+1}}  h_{\sigma_{n+1}-}$,
so that
$$
\xi_{\sigma_{n+1}-}=\xi_{\sigma_{n}}+\int_{\sigma_{n}}^{\sigma_{n+1}}
(u_{s}[h\hat h],L_{s}\zeta(h_{s}\cdot))\,ds
$$
$$
+\int_{\sigma_{n}}^{\sigma_{n+1}}(h_{s}\hat h_{s}f_{s},\zeta(h_{s}\cdot))\,ds,
$$
$$
\xi_{\sigma_{n+1} }=\xi_{\sigma_{n+1}-}+
[(u_{\sigma_{n+1}}[h\hat h],\zeta(g_{\sigma_{n+1}}   
h_{\sigma_{n+1}-}\cdot))
 -\xi_{\sigma_{n+1}-}].
$$
It follows easily that  for  (each $\omega$ and)  $t\in [0,T]$ we have
$$
(u_{t}[h\hat h],\zeta(h_{t}\cdot)) 
=(u_{0},\zeta) 
+\int_{0}^{t}
(u_{s}[h\hat h],L_{s}\zeta(h_{s}\cdot)) \,ds
$$
$$
+\int_{0}^{t}(h_{s}\hat h_{s}f_{s},\zeta(h_{s}\cdot))
\,ds
$$
\begin{equation}
                                                  \label{5.15.3}
+\int_{(0,t]}\big [(u_{s}[h\hat h],\zeta(g_{s}h_{s-}\cdot))
-\xi_{s-}\big]\,d\pi_{s}.
\end{equation}
The above formulas show that $\xi_{t-}$ is a well-defined left-continuous
process, which is $\cF_{t}$-adapted since $\xi_{t}$ is such
(cf. Lemma \ref{lemma 2.28.3}).
We observe also that, by Remark \ref{remark 5.15.2}
and the fact that $h_{t-}$ is left-continuous, $\cF_{t}$-adapted,
and hence predictable process, the last integrand is
predictable. 

Of course, we want to take expectations of both sides of \eqref{5.15.3}
and use Lemma \ref{lemma 2.29.2}. 
Introduce
$$
v_{t}(x)=u_{t}[h(\omega)\hat h]
(h^{-1}_{t }(\omega) x).
$$
By Lemma \ref{lemma 2.28.3} we have $v\in W$  
for any $\omega$. In particular, 
 $|v_{t}(x)|\leq K$
   Hence, changing variables 
 (see Remark \ref{remark 8.24.1}) we find
$$
E\int_{0}^{t}|\xi_{s-}|\,ds=
E\int_{0}^{t}|\xi_{s}|\,ds=
E\int_{0}^{t}|(v_{s},\zeta)|\,ds
\leq  K T \int_{\bR^d} |\zeta(y)| dy  
<\infty.
$$
Similarly,  
$$
E\int_{0}^{t}( |u_{s}[h\hat h] |,
|\zeta(g_{s}h_{s-}\cdot) |)
|\,ds  \le   E\int_{[0,T]
\times\bR^{d}}|v_{s}(g^{-1}_{s}h^{-1}_{s-}x)| |\zeta(x)|\,dxds   
$$
$$
\leq KT \int_{\bR^d} |\zeta(y)| dy   
<\infty.
$$
Dealing with other terms on the right in \eqref{5.15.3}
presents no problem either, and, after taking the expectations
of both sides and using Fubini's theorem, we obtain
$$
E(u_{t}[h\hat h](h^{-1}_{t}\cdot),\zeta)=
(u_{0},\zeta)+\int_{0}^{t}(Ev_{s},L_{s}\zeta)\,ds
+\int_{0}^{t}(\hat h_{s}f_{s},\zeta)\,ds
$$
\begin{equation}
                                                          \label{5.16.1}
+\lambda\int_{0}^{t} (Ev_{s}(g^{-1}_{s}\cdot)-Ev_{s},\zeta)\,ds.
\end{equation}

Since $P(t\in\{\sigma_{1},\sigma_{2},...\})=0$, we have
$E(u_{t}[h\hat h](h^{-1}_{t}\cdot),\zeta)=E(v_{t} ,\zeta)$.
Furthermore,
$$
E(|v_{t}| ,|\zeta|)\leq  K  \int_{\bR^d} |\zeta(y)|dy 
 <\infty,
$$
which allows us to use Fubini's theorem and conclude
that $E(v_{t} ,\zeta)=(Ev_{t} ,\zeta)$. This and \eqref{5.16.1}
show that the function 
$$
w_{t}[\hat h](x)=w_{t}(x):=Ev_{t}(x)
=Eu_{t}[h \hat h](h^{-1}_{t }  x),
$$
 which belongs to $W$ by Lemma \ref{lemma 2.28.3} and by assumption,
satisfies
$$ w_{t}(x)=u_{0}(x)+
\int_{0}^{  t}[L^{*}_{r}w_{r}(x)+\lambda(g^{-1}_{r}-1)w_{r}(x)
+\hat h_{r}f_{r}(x)]\,dr.
$$
This equation coincides with  
$$
w_{t}(x)=u_{0}(x)+
\int_{0}^{ t}[\hat L^{*}_{r}w_{r}(x)
+\hat h_{r}f_{r}(x)]\,dr
$$ 
if we take ($n=1$ and) $g^{(1)}_{t}=  g^{-1}_{t}  $,
  which is as arbitrary
as a member of $B((0,T),\cG)$ could be. Hence
Assumption \ref{assumption 2.27.1} (iii) is satisfied. Finally,
if $h',h''\in B((0,T),\cG)$, then due to our assumptions
$$
|w_{t}[ h'](x)-w_{t}[  h''](x)| 
$$
$$
=|Eu_{t}[h h' ](h^{-1}_{t }  x)-Eu_{t}[h   h'' ](h^{-1}_{t }  x)|
\leq
\sup_{\omega,x}
|u_{t}[h h' ](   x)- u_{t}[h   h'' ](  x)|
$$
$$
\leq   K \sup_{\omega,x}\int_{0}^{t}
\sup_{y\in\bR^{d}}|f_{r}(  h'_{r}h_{r} y)- 
f_{r}(h''_{r}h_{r} y)|\,dr =
  K   \int_{0}^{t}
\sup_{y\in\bR^{d}}|f_{r}(  h'_{r}  y)- 
f_{r}(h''_{r} y)|\,dr ,
$$
which shows that Assumption \ref{assumption 2.27.1} (iv) is
satisfied as well and proves the theorem.

\mysection{Proof of Theorems 
\protect\ref{theorem 9.17.1} and \protect\ref{theorem 2.29.1}}

                                       \label{section 6.29.6}

{\bf Proof of Theorem \ref{theorem 9.17.1}}.
 By Assumption \ref{assumption 2.26.3} (i)
all $u^{k}_{t}(x)$, $k\geq1$, are uniformly bounded on
$[0,T]\times\bR^{d}$. Then there exists a subsequence
still denoted by $u^k$ and a bounded (Borel) function
$u$ on $[0,T]\times\bR^{d}$ such that for any
$\zeta\in L_{1}( [0,T]\times\bR^{d})$ we have
$$
\int_{ [0,T]\times\bR^{d} }u^{k}_{t}(x)\zeta_{t}(x)\,dxdt
\to \int_{ [0,T]\times\bR^{d} }u _{t}(x)\zeta_{t}(x)\,dxdt.
$$
 Next, take $\zeta\in C^{\infty}_{0}(\bR^{d})$ and write
that by definition
\begin{equation}
                                                \label{9.17.2}
 (u^{k}_{t},\zeta) =(u_{0},\zeta)
+\int_{0}^{t}\big[(u^{k}_{s},L^{0}_{s}\zeta)
 +
( f_{s},\zeta)\big]\,ds+F^{k}_{t},
\end{equation}
where
$$
F^{k}_{t}=\int_{0}^{t}\big (u^{k}_{s},(L^{k}_{s}-L^{0}_{s})
\zeta\big)\,ds.
$$
  Let us fix $t \in (0,T]$.  In light of Assumption \ref{assumption 2.26.3} (i)
and requirement b) in the theorem we have that
$F^{n}_{t}\to0$ as $n\to\infty$. 

There are two consequences of this fact. First,
the right-hand sides of \eqref{9.17.2} converge
as $k \to\infty$ to 
$$
(u_{0},\zeta)
+\int_{0}^{t}\big[(u _{s},L^{0}_{s}\zeta)
 +
( f_{s},\zeta)\big]\,ds.
$$
Secondly, the  left-hand  sides of
\eqref{9.17.2}  
also converge 
  for any
$\zeta\in C^{\infty}_{0}(\bR^{d})$ 
 to the limit, say $\phi_{t}(\zeta)$,
which is a generalized function.
  Since
$$
|\phi_{t}(\zeta)|\leq K\int_{\bR^{d}}|\zeta   (x)  |\,dx,
$$ 
  $\phi_t$ can be extended to a 
linear continuous functional on $L^1(\R^d)$ and so  
there exists
 a (bounded) function
$u=u_{t}(x)$, $(t,x)\in[0,T]\times\bR^{d}$, such that
\begin{equation}
                                          \label{8.27.2}
\phi_{t}(\zeta)=
\int_{\bR^{d}}u _{t}(x)\zeta(x)\,dx,\quad
\int_{\bR^{d}}u^{k}_{t}(x)\zeta(x)\,dx\to 
\int_{\bR^{d}}u _{t}(x)\zeta(x)\,dx
\end{equation}
 for any  
  $t\in[0,T]$ and $\zeta\in C^{\infty}_{0}(\bR^{d})$. 

 Another way to get the same result
is to fix $R>0$, take the ball $B$
of radius $R$ centered at the origin,
and take a subsequence $u^{k'}_{t}$ such
that $u^{k'}_{t}I_{B}$   converges  weakly in 
$L_{2}(B)$ to a function $u^{B}_{t}$. Then,
obviously, 
\begin{equation}
                                          \label{8.27.1}
\int_{\bR^{d}}u^{B}_{t}  (x)  \zeta   (x) \,dx= \phi_t (\zeta)
\end{equation}
for any $\zeta\in C^{\infty}_{0}(B)$. This holds
for any weakly convergent subsequence
of $u^{k}_{t}I_{B}$, and shows that the weak limit
is always the same. Hence, the whole
sequence $u^{k}_{t}I_{B}$ converges weakly
in $L_{2}(B)$ to $u^{B}_{t}$. Of course, 
\eqref{8.27.1} implies that, for balls $B'\subset B''$,
$u_{t}^{B'}=u_{t}^{B''}$ on $B'$ and this allows us
to define $u_{t}$ on $\bR^{d}$ for which 
\eqref{8.27.2} hold  for any $\zeta\in C^{\infty}_{0}
(\bR^{d})$.

  By Assumption \ref{assumption 5.17.1}
there exists $w\in W$ such that
$w_{t}=u_{t}$ \(a.e.\) on $\bR^{d}$ for any $t\in[0,T]$. 
It follows that 
\begin{equation}
                                          \label{8.24.4}
\int_{\bR^{d}}u^{k}_{t}(x)\zeta(x)\,dx\to 
\int_{\bR^{d}}w_{t}(x)\zeta(x)\,dx
\end{equation}
 for any  
  $t\in[0,T]$ and $\zeta\in C^{\infty}_{0}(\bR^{d})$. 
Hence, for any $t\in[0,T]$ and $\zeta\in C^{\infty}_{0}(\bR^{d})$
\begin{equation}
                                          \label{8.24.3}
(w_{t},\zeta)=(u_{0},\zeta)+  \int_{0}^{t}\big[(u_{s},
L^{0}_{s}\zeta 
  )+
(h_{s}f_{s},\zeta)\big]\,ds.
\end{equation}

Next, note that, for any smooth function $\eta_{t}(x)$
with compact support in $(0,T)\times\bR^{d}$, on the one hand,
by definition of $u_{t}(x)$
$$
\lim_{k\to\infty}
\int_{(0,T)\times\bR^{d}}u^{k}_{t}(x)\eta_{t}(x)\,dxdt
=\int_{(0,T)\times\bR^{d}}u _{t}(x)\eta_{t}(x)\,dxdt.
$$
On the other hand, owing to \eqref{8.24.4},
by the dominated convergence theorem,
$$
\lim_{k\to\infty}
\int_{0}^{T}\,dt\int_{\bR^{d}}u^{k}_{t}(x)\eta_{t}(x)\,dx
=\int_{0}^{T}\,dt\int_{\bR^{d}}w_{t}(x)\eta_{t}(x)\,dx.
$$
It follows that $u_{t}(x)=w_{t}(x)$ (a.e.) in 
$(0,T)\times\bR^{d}$ and we can replace $u_{s}$
with $w_{s}$ in \eqref{8.24.3} without violating this equality.
This proves the theorem.    \qed

We build our proof
of Theorem \ref{theorem 2.29.1}
 entirely on Theorems \ref{theorem 2.27.1}
 and \ref{theorem 9.17.1}
thus avoiding using probability theory. 
We need the following. 
Recall that if $\nu\in\frN$, we set $\nu^{0}x=\nu x-\nu0$,
and $\nu^{0}$ is a linear mapping.

\begin{lemma}
                                              \label{lemma 5.17.1}
Let $\nu\in\frN$ and let $g$ be the solution of
\eqref{2.27.5}. Then,
for any $\zeta\in C^{\infty}_{0}(\bR^{d})$ 
and $x\in\bR^{d}$, 
$$
\zeta(g^{-1}(t)x) =\zeta(x)- (\nu x)^{i}D_{i}\zeta(x) t 
+(1/2)\big[(\nu x)^{i}(\nu x)^{j}D_{ij}\zeta(x)
$$
$$
+(\nu^{0}\nu x)^{i} D_{i}\zeta(x) 
 \big]t^{2}+o(t^{2})
$$ 
as $t\downarrow0$.
\end{lemma}

Proof. For any $y\in\bR^{d}$ we have $\dot g(t)y=\nu^{0} g(t)y+\nu0$.
The solution of this equation which equals $y$ at $t=0$ is
$$
g(t)y=e^{\nu^{0}t}y+\int_{0}^{t}e^{\nu^{0}s}\nu0\,ds.
$$
It follows that
$$
  g^{-1}(t) x=e^{-\nu^{0}t}x-\int_{0}^{t}e^{\nu^{0}(s-t)}\nu 0\,ds
$$ 
and the results follow by Taylor's formula. 
The lemma is proved.     \qed

{\bf Proof of Theorem \ref{theorem 2.29.1}}.
Take $\mu^{(1)},...,  \mu^{(n)  } \in B((0,T),\frN)$
and for $k=1,2,...$, set
$$
L^{k}_{t}=L_{t}+\sum_{r=1}^{n}M^{(r)k}_{t}.
$$
where
$$
M^{(r)k}_{s}\phi(x)
:=k^{2}[\phi (e^{\mu^{(r)}_{s}/k}x)-2\phi (x)
+\phi (e^{-\mu^{(r)}_{s}/k}x)].
$$
Observe that $M^{(r)k}_{s}$ are formally self-adjoint,
so that
$$
L^{k*}_{t}=L^{*}_{t}+\sum_{r=1}^{n}M^{(r)k}_{t}
$$
and by Theorem \ref{theorem 2.27.1}, for any $k\geq 1$,
and $h\in B((0,T),\cG)$ there exists $u^{k}\in W$
satisfying 
 \begin{equation} 
                                      \label{9.17.3}  
(u^{k}_{t},\zeta)  
=(u_{0},\zeta)+\int_{0}^{t}(u^{k}_{s}, L^{k}_{s} \zeta)
\,ds
+\int_{0}^{t}  ( h_{s}f_{s}, \zeta) \,ds.
\end{equation}
for any $\zeta\in C^{\infty}_{0}(\bR^{d})$ and $t\in[0,T]$.

Then define
$$
L^{0}_{t}=L_{t}+\sum_{r=1}^{n} M_{\mu^{(r)}_{t}}  .
$$
Observe that owing to the boundedness of the $\mu^{(r)}$'s,
it follows easily from the arguments in the proof of Lemma
\ref{lemma 5.17.1} that there is a ball $B$ such that
$M^{(r)k }_{s}\zeta=0$ outside  $B$ for all $k$ 
and $s\in(0,T)$ and   
$$
M^{(r)k }_{s}\zeta(x)\to
  M_{\mu^{(r)}_{s}}\zeta(x)   
$$
as $k\to\infty$
 uniformly with respect to $s\in(0,T)$ and $x\in B$.
It now follows by Theorem \ref{theorem 9.17.1}
that there exists $u^{0}\in W$ for which
\eqref{9.17.3} holds with   $ 0$
 in place of   $ k$  
for any $\zeta\in C^{\infty}_{0}(\bR^{d})$ and $t\in[0,T]$.
 This is exactly what we need because  
simple manipulations show that, for $\nu\in \frN$,
$$
M_{\nu}\zeta(x)=(\nu x)^{i}D_{i}\big[(\nu x)^{j}D_{j}\zeta(x)
],
$$
so that the operators
$M_{\nu}$ are formally self-adjoint, 
and this proves the theorem. \qed

\section{Possible extensions  
to non-local operators}
                                 \label{section 9.17.1}

\begin{assumption}
                                \label{assumption 9.17.1}

We are given a family $\{\nu_{t}(A),t\in(0,T)\}$
of   measures on Borel subsets of $\bR^{d} $
such that 

(i) $\nu_{t}(\{0\})=0$ for any $t\in(0,T)$,

(ii) $\nu_{t}(A)$ is a (Borel) measurable function
of $t\in(0,T)$,

(iii) we have
$$
\int_{\bR^{d}}(1\wedge|x|^{2})\,\nu_{t}(dx)<\infty\quad\forall
t\in(0,T),\quad
\int_{(0,T)\times\bR^{d}}(1\wedge|x|^{2})\,\nu_{t}(dx)dt
<\infty.
$$
\end{assumption} 
\begin{assumption}
                                \label{assumption 9.17.2}

We are given  
$W$, $\mathcal G$, $K$, $L$, $  u_0 $,  $f$
 as in Theorem \ref{theorem 2.29.1} with $\mathcal G$ being the group
of translations. 
\end{assumption} 

Introduce
\begin{equation}
                                          \label{9.17.4}
L^{0}_{t}=L_{t}+J_{\nu_{t}},
\end{equation}
where, for $\phi\in C^{\infty}_{0}(\bR^{d})$ and measure $\nu$,
$$
J_{\nu }\phi(x)=\int_{\R^d} \Big[\phi(x+y) - \phi (x) - 
y^{i} D_{i}\phi ( x) \, I_{\{ |y| \le 1 \}} (y)  \Big]  \nu  (dy).
$$
As a side observation
 recall that if $\nu$ is a measure on $\bR^{d}\setminus
\{0\}$
such that
$$
\int_{\bR^{d}}(1\wedge|x|^{2})\,\nu (dx)<\infty,
$$
the operator $J_{\nu}$
is known in probability theory as
 the generator of a unique in law L\'evy process
associated to $\nu$  (this process is without Gaussian part; 
see \cite{sato} and \cite{Kr02}).

One knows (and we will see this again in the proof of
Theorem \ref{theorem 9.17.2}) that,
owing to Assumption \ref{assumption 9.17.1}, $J_{\nu_{t} }\phi(x)$
is well defined for any $\phi\in C^{\infty}_{0}(\bR^{d})$.
Standard measure theoretic arguments show that $J_{\nu_{t} }\phi(x)$
is a measurable function of $(t,x)$ for any
$\phi\in C^{\infty}_{0}(\bR^{d})$.

\begin{theorem}
                                            \label{theorem 9.17.2}
Under the above assumptions for any $h\in B((0,T),\cG)$
there exists $u\in W$ such that
\eqref{weak11} holds for any $\zeta\in C^{\infty}_{0}(\bR^{d})$ and $t\in[0,T]$
with $L^{0}_{s}$ in place of $L_{s}$.  
\end{theorem}

Proof. Notice that by Taylor's formula for any
$\phi\in C^{\infty}_{0}(\bR^{d})$, if $|y|\leq1$, then
$$
|\phi(x+y) - \phi (x) - y^{i} D_{i}\phi ( x) |\leq |y|^{2}N\sup_{|z|\leq1}|D^{2}\phi(x+z)|,
$$
where $N$ is a constant.
  Below by $N$ we denote generic constants
which may change from one occurrence to another. 
It follows that
$$
|J_{\nu_{t} }\phi(x)|\leq N\sup_{|z|\leq1} |D^{2}\phi(x+z)|
\int_{|y|\leq1} |y|^{2} \,\nu_{t}(dy)
$$
$$
+\int_{|y|\geq1}(| \phi(x+y)|+|\phi(x)|) \,\nu_{t}(dy),
$$
and owing to Assumption \ref{assumption 9.17.1} (iii)
and Fubini's theorem
we see that Assumption \ref{assumption 2.27.1} (ii) is satisfied
with $L^{0}_{t}$ in place of $L_{t}$.

Furthermore, for Borel sets $A\subset\bR^{d}$,
 define $\nu^{k}_{t}(A)=\nu_{t}(A\cap B_{k})$, 
where $B_{k}=\{|x|\leq k\}$. Then the above manipulations show that,
for $k\geq1$
$$
\delta^{k}:=\int_{ (0,T)\times\bR^{d} } |(J_{\nu_{t} }-J_{\nu^{k}_{t} })\phi(x)| 
\,dt dx
$$
$$
\leq \int_{ (0,T)\times\bR^{d} }
\int_{|y|\geq k}(| \phi(x+y)|+|\phi(x)|) \,\nu_{t}(dy)\,dt dx
$$
$$
\leq N\int_{(0,T)}\nu_{t}(B_{k}^{c})\,dt,
$$
which tends to zero as $k\to\infty$ by the dominated convergence theorem
 (see (iii) in Assumption \ref{assumption 9.17.1}). 
It follows that, if we introduce $L^{k}_{t}$
by \eqref{9.17.4} with $J_{\nu^{k}_{t}}$ in place of
$J_{\nu _{t}}$, then condition b) of Theorem
 \ref{theorem 9.17.1} is fulfilled. Of course, 
condition a) is fulfilled as well by the above.
Now thanks to Theorem
 \ref{theorem 9.17.1} to prove our theorem, it suffices
to prove it with $\nu^{k}_{t}$ in place of $\nu_{t}$.

Hence, below we assume that $\nu_{t}(B^{c})=0$, where $B$ is
a ball (independent of $t$). We can play the same trick
with small jumps. Set this time 
$\nu^{k}_{t}(A)=\nu_{t}(A\cap B^{c}_{1/k})$
(of course, this $\nu^{k}_{t}$ is different from the above one,
but it is convenient to forget the above
$\nu^{k}_{t}$ and introduce $\delta^{k}$ by the same formula
with the new $\nu^{k}_{t}$). Then
$$
\delta^{k}\leq N\int_{ [0,T]\times\bR^{d} }
\sup_{|z|\leq1} |D^{2}\phi(x+z)|
\int_{|y|\leq1/k} |y|^{2} \,\nu_{t}(dy)\,dt dx
$$
$$
\leq  N \int_{(0,T)}\int_{|y|\leq1/k} |y|^{2} \,\nu_{t}(dy)\,dt,
$$
which again tends to zero as $k\to\infty$ 
 by the dominated convergence theorem.

Since this measures $\nu^{k}_{t}$ are finite,
we now see that we may concentrate on the case in which
$\nu_{t}$ are finite measures with support in a ball $B$ independent of
$t$. One more simplification is   achieved 
by introducing 
$$
\nu^{k}_{t}=\nu_{t}I_{ \{  \nu_{t}(\bR^{d})\leq k \} },
$$
in which case
$$
\delta^{k}=
\int_{ (0,T)\times\bR^{d} } | J_{\nu_{t} } \phi(x)| 
I_{ \{  \nu_{t}(\bR^{d})\geq k  \}}
\,dt dx\to0
$$
by the dominated convergence theorem.

Thus, we need only consider the case in which
$  t \mapsto  \nu _{t}(\bR^{d})$ is bounded and 
$\nu_{t}$ have support in a ball $ B = \{|x|\leq R\}$.

In that case
$$
J_{\nu_{t} }\phi(x)=\int_{B}[\phi(x+y)-\phi(x)]\,\nu_{t}(dy)
+b_{t}\cdot D\phi(x),
$$
where
$$
b_{t}=\int_{|y|\leq1}y\,\nu_{t}(dy).
$$
 For
$y\in\bR^{d}$ set $\kappa_{k}(y)=(\kappa_{k}(y^{1}),...,
\kappa_{k}(y^{d}))$ and introduce
$$
J^{k}_{\nu_{t} }\phi(x)=
\int_{B}[\phi(x+\kappa_{k}(y))-\phi(x)]\,\nu_{t}(dy)
+  k[\phi(x+b_{t}/k)-\phi(x)].
$$
As is (very) easy to see 
$$
 \delta^{k} = 
\int_{ (0,T)\times\bR^{d} } |(J_{\nu_{t} }-J_{\nu_{t} }^k)\phi(x)| 
\,dt dx
\to 0  
$$ 
as $k\to\infty$
and to finish the proof it only remains to
refer to Theorem \ref{theorem 10.6.2} after observing 
  that
$$
J^{k}_{\nu_{t} }\phi(x)=\sum_{\substack{z\in(1/k)\bZ^{d}\\
|z|\leq R+1}}[\phi(x
+z)-\phi(x)]\nu_{t}(\{y:\kappa_{k}(y)=z\})
+  k[\phi(x+b_{t}/k)-\phi(x)],
$$
where the sum contains only finite number of  terms.
The theorem is proved.    \qed

   \textbf{Acknowledgements.}
  The work on this  paper   started during the conference 
   ``New advances in PDE's, Inverse Problems and Control Theory'' 
 (July 6-10, 2015 Parma). The authors would like to thank the  Mathematical 
Department of Parma.

\end{document}